\documentclass[10pt]{amsart}
\usepackage[all]{xy}
\usepackage{xspace}
\usepackage{bm}
\usepackage{amsmath,amsthm}
\usepackage{amstext}
\usepackage{amsfonts}
\usepackage[mathscr]{euscript}
\usepackage{amscd}
\usepackage{amssymb}
\usepackage{euscript}
\begin{document}
%
%
\theoremstyle{plain}

\swapnumbers
    \newtheorem{thm}{Theorem}[section]
    \newtheorem{prop}[thm]{Proposition}
    \newtheorem{lemma}[thm]{Lemma}
    \newtheorem{cor}[thm]{Corollary}
    \newtheorem{fact}[thm]{Fact}
        \newtheorem*{thma}{Theorem A}
        \newtheorem*{thmb}{Theorem B}
        \newtheorem*{thmc}{Theorem C}
    \newtheorem{subsec}[thm]{}
\theoremstyle{definition}
    \newtheorem{assume}[thm]{Assumption}
    \newtheorem{defn}[thm]{Definition}
    \newtheorem{example}[thm]{Example}
    \newtheorem{examples}[thm]{Examples}
    \newtheorem{claim}[thm]{Claim}
    \newtheorem{notn}[thm]{Notation}
    \newtheorem{summary}[thm]{Summary}
    \newtheorem{construct}[thm]{Construction}
\theoremstyle{remark}
\newtheorem{remark}[thm]{Remark}
        \newtheorem{aside}[thm]{Aside}
    \newtheorem{ack}[thm]{Acknowledgments}
\newenvironment{myeq}[1][]
{\stepcounter{thm}\begin{equation}\tag{\thethm}{#1}}
{\end{equation}}
\newcommand{\mydiag}[2][]{\myeq[#1]{\xymatrix{#2}}}
\newcommand{\mydiagram}[2][]
{\stepcounter{thm}\begin{equation}
     \tag{\thethm}{#1}\vcenter{\xymatrix{#2}}\end{equation}}
\newcommand{\mydiagr}[2][]
{\stepcounter{thm}\begin{equation}
     \tag{\thethm}{#1}\vcenter{\xymatrix@R=22pt{#2}}\end{equation}}
\newcommand{\mydiags}[2][]
{\stepcounter{thm}\begin{equation}
     \tag{\thethm}{#1}\vcenter{\xymatrix@R=17pt{#2}}\end{equation}}
%
\newcommand{\mypict}[2][]
{\stepcounter{thm}\begin{equation}
     \tag{\thethm}{#1}\vcenter{\begin{picture}{#2}\end{picture}}\end{equation}}
\newenvironment{mysubsection}[2][]
{\begin{subsec}\begin{upshape}\begin{bfseries}{#2.}
\end{bfseries}{#1}}
{\end{upshape}\end{subsec}}
\newenvironment{mysubsect}[2][]
{\begin{subsec}\begin{upshape}\begin{bfseries}{#2\vsm.}
\end{bfseries}{#1}}
{\end{upshape}\end{subsec}}
\newcommand{\supsect}[2]
{\vspace*{-5mm}\quad\\ \begin{center}\textbf{{#1}}\vsm.~~~~\textbf{{#2}}\end{center}}
\newcommand{\sect}{\setcounter{thm}{0}\section}
\newcommand{\wh}{\ -- \ }
\newcommand{\wwh}{-- \ }
\newcommand{\w}[2][ ]{\ \ensuremath{#2}{#1}\ }
\newcommand{\ww}[1]{\ \ensuremath{#1}}
\newcommand{\www}[2][ ]{\ensuremath{#2}{#1}\ }
\newcommand{\wb}[2][ ]{\ (\ensuremath{#2}){#1}\ }
\newcommand{\wwb}[1]{\ (\ensuremath{#1})-}
\newcommand{\wref}[2][ ]{\ \eqref{#2}{#1}\ }
%
%
\newcommand{\xra}[1]{\xrightarrow{#1}}
\newcommand{\xla}[1]{\xleftarrow{#1}}
\newcommand{\lar}{\leftarrow}
\newcommand{\Ra}{\Rightarrow}
\newcommand{\hra}{\hookrightarrow}
\newcommand{\lora}{\longrightarrow}
\newcommand{\bstar}{\mbox{\large $\star$}}
\newcommand{\adj}[2]{\substack{{#1}\\ \rightleftharpoons \\ {#2}}}
\newcommand{\hsq}{\hspace*{13 mm}}
\newcommand{\hsp}{\hspace*{10 mm}}
\newcommand{\hs}{\hspace*{5 mm}}
\newcommand{\hsm}{\hspace*{2 mm}}
\newcommand{\vsp}{\vspace{10 mm}}
\newcommand{\vs}{\vspace{7 mm}}
\newcommand{\vsm}{\vspace{2 mm}}
\newcommand{\rest}[1]{\lvert\sb{#1}}
\newcommand{\lo}[1]{\sb{(#1)}}
\newcommand{\up}[1]{\sp{(#1)}}
\newcommand{\bp}[1]{\sp{[#1]}}
\newcommand{\lp}[1]{\sb{[#1]}}
\newcommand{\lra}[1]{\langle{#1}\rangle}
\newcommand{\lin}[1]{\{{#1}\}}
\newcommand{\EQUIV}{\Leftrightarrow}
\newcommand{\epic}{\lora\hspace{-3 mm}\to}
\newcommand{\xepic}[1]{\stackrel{#1}{\lora}\hspace{-3 mm}\to}
%
%
\newcommand{\q}[1]{^{({#1})}}
\newcommand{\hp}[1]{\widehat{p}\q{#1}}
\newcommand{\hj}[1]{\widehat{j}\q{#1}}
\newcommand{\li}[1]{\sb{({#1})}}
\newcommand{\od}{\overline{\partial}}
\newcommand{\ud}{\overline{d}}
\newcommand{\odz}[2]{\partial\sp{#1}\sb{#2}}
\newcommand{\udz}[1]{\ud\sp{0}\sb{#1}}

\newcommand{\wrh}{\hat{\rho}}
\newcommand{\vare}{\varepsilon}
\newcommand{\vpsi}[1]{\psi\sb{#1}}
\newcommand{\tpsi}[1]{\widetilde{\psi}\sb{#1}}
\newcommand{\varp}[1]{\varphi\sb{#1}}
\newcommand{\tvarp}[1]{\widetilde{\varphi}\sb{#1}}
\newcommand{\vart}{\vartheta}
%
%
\newcommand{\ab}{\operatorname{ab}}
\newcommand{\Ab}{\operatorname{Ab}}
\newcommand{\Aut}{\operatorname{Aut}}
\newcommand{\cf}{\operatorname{cof}}
\newcommand{\Cok}{\operatorname{Coker}\,}
\newcommand{\colim}{\operatorname{colim}}
\newcommand{\colimit}[1]{\raisebox{-1.7ex}{$\stackrel{\colim}{\scriptstyle{#1}}$}}
\newcommand{\csk}[1]{\operatorname{csk}\sb{#1}}
\newcommand{\Dec}{\operatorname{Dec}}
\newcommand{\diag}{\operatorname{diag}}
\newcommand{\ev}{\operatorname{ev}}
\newcommand{\Ext}{\operatorname{Ext}}
\newcommand{\fib}{\operatorname{fib}}
\newcommand{\fold}{\nabla}
\newcommand{\free}{\sp{\operatorname{free}}}
\newcommand{\gr}{\operatorname{gr}}
\newcommand{\ho}{\operatorname{ho}}
\newcommand{\Hom}{\operatorname{Hom}}
\newcommand{\Id}{\operatorname{Id}}
\newcommand{\Image}{\operatorname{Im}\,}
\newcommand{\inc}{\operatorname{inc}}
\newcommand{\Ind}{\operatorname{Ind}}
\newcommand{\Kan}{\operatorname{Kan}}
\newcommand{\Ker}{\operatorname{Ker}\,}
\newcommand{\map}{\operatorname{map}}
\newcommand{\mapa}{\map\sb{\ast}}
\newcommand{\Obj}{\operatorname{Obj}}
\newcommand{\op}{\sp{\operatorname{op}}}
\newcommand{\proj}{\operatorname{proj}}
\newcommand{\real}{\operatorname{re}}
\newcommand{\red}{\operatorname{red}}
\newcommand{\rel}{\operatorname{rel}}
\newcommand{\res}{\operatorname{res}}
\newcommand{\sk}[1]{\operatorname{sk}\sb{#1}}
\newcommand{\spec}{\operatorname{sp}}
\newcommand{\Stov}{\operatorname{St}}
\newcommand{\Tor}{\operatorname{Tor}}
\newcommand{\Tot}{\operatorname{Tot}}
\newcommand{\we}{\operatorname{w.e.}}
\newcommand{\sss}{\hspace*{1 mm}\sp{s}}
\newcommand{\vvv}{\hspace*{1 mm}\sp{v}}
%
%
\newcommand{\A}{{\EuScript A}}
\newcommand{\hA}{\hat{\A}}
\newcommand{\uA}{{\EuScript A}}
\newcommand{\whA}{\widehat{\A}}
\newcommand{\tA}{\tilde{A}}
\newcommand{\cA}{{\mathcal A}}
\newcommand{\AR}{\A\sb{R}}
\newcommand{\B}{{\EuScript B}}
\newcommand{\cB}{{\mathcal B}}
\newcommand{\hB}{\hat{\B}}
\newcommand{\AB}{(\A,\B)}
\newcommand{\tB}{\tilde{\B}}
\newcommand{\C}{{\mathcal C}}
\newcommand{\Cab}{\C\sb{\ab}}
\newcommand{\hC}{\hat{\C}}
\newcommand{\tC}{\tilde{\C}}
\newcommand{\Tow}[1]{\hspace*{-0.5mm}\sp{\#}{#1}}
\newcommand{\Ct}{\Tow{\C}}
\newcommand{\Ctl}[1]{\hspace*{-1mm}\sb{{#1}\leq}\sp{\hspace*{1.7mm}\#}
                 \hspace{0.5mm}\C}
\newcommand{\Ctg}[1]{\hspace*{-1mm}\sb{{#1}\geq}\sp{\hspace*{1.7mm}\#}
                  \hspace{0.5mm}\C}
\newcommand{\D}{{\mathcal D}}
\newcommand{\hD}{\hat{\D}}
\newcommand{\E}{{\mathcal E}}
\newcommand{\Enr}[2]{\Stem\sp{#1}\sb{{#2}}}
\newcommand{\Ens}[1]{\Enr{#1}{\bDel\op}}
\newcommand{\Enc}[1]{\Enr{#1}{\bDel}}
\newcommand{\F}{{\mathcal F}}
\newcommand{\G}{{\mathcal G}}
\newcommand{\Gn}[1]{\G\sb{[#1]}}
\newcommand{\K}{{\mathcal K}}
\newcommand{\LL}{{\mathcal L}}
\newcommand{\LR}{\LL\sp{R}}
\newcommand{\LS}{\LL\sb{\uS}}
\newcommand{\M}{{\mathcal M}}
\newcommand{\sMR}{s\M\sp{R}}
\newcommand{\sMRl}[1]{\sMR\sb{#1}}
\newcommand{\sMRla}{\sMRl{\lambda}}
\newcommand{\MA}{\M\sb{\A}}
\newcommand{\MAd}{\MA\sp{\delta}}
\newcommand{\hM}{\hat{\M}}
\newcommand{\hMA}{\hM\sb{\A}}
\newcommand{\hMAp}{\hM\sb{\A'}}
\newcommand{\tM}{\tilde{\M}}
\newcommand{\tMA}{\tM\sb{\A}}
\newcommand{\nMA}[1]{\hy{#1}{\M\sb{\A}}}
\newcommand{\tMuA}{\tM\sp{\A}}
\newcommand{\N}{{\mathcal N}}
\newcommand{\OO}{{\mathcal O}}
\newcommand{\PP}{{\mathcal P}}
\newcommand{\cQ}{{\mathcal Q}}
\newcommand{\Qo}[1]{\cQ[#1]}
\newcommand{\cQd}{\cQ\sb{\bullet}}
\newcommand{\hQd}{\widehat{\cQ}\sb{\bullet}}
\newcommand{\R}{{\mathcal R}}
\newcommand{\Rd}{\R\sb{\bullet}}
\newcommand{\hRd}{\widehat{\R}\sb{\bullet}}
\newcommand{\RR}{\R\sp{R}}
\newcommand{\RS}{\R\sb{\uS}}
\newcommand{\uS}{{\EuScript B}}
\newcommand{\eS}{{\EuScript S}}
\newcommand{\hS}[1]{\hat{\eS}\sp{#1}}
\newcommand{\Ss}{{\mathcal S}}
\newcommand{\wS}{\widehat{S}}
\newcommand{\cS}{{\mathcal S}}
\newcommand{\Sd}{\cS\sb{\bullet}}
\newcommand{\TT}{{\mathcal T}}
\newcommand{\TS}{\TT\sb{\uS}}
\newcommand{\TR}{\TT\sp{R}}
\newcommand{\U}{{\mathcal U}}
\newcommand{\V}{{\mathcal V}}
\newcommand{\VS}{\V\sb{\uS}}
\newcommand{\VR}{\V\sp{R}}
\newcommand{\X}{{\mathcal X}}
\newcommand{\cX}{{\mathcal X}}
\newcommand{\Wc}{\cW'}
\newcommand{\Wf}{\cW\sb{\F}}
%
%
\newcommand{\Stem}{\mbox{\sf Stem}}
\newcommand{\Sys}{\mbox{\sf Sys}}
\newcommand{\Pstem}{\mbox{\sf PStem}}
\newcommand{\Po}[1]{P\sp{#1}}
\newcommand{\PPo}[1]{\PP[{#1}]}
\newcommand{\Pnk}[3]{\Po{#1}\wk{#2}{#3}}
\newcommand{\Ql}[1]{Q\lo{#1}}
\newcommand{\Qu}[1]{Q\up{#1}}
%
%
\newcommand{\hy}[2]{{#1}\text{-}{#2}}
\newcommand{\Al}{\mbox{\sf Alg}}
\newcommand{\Alg}[1]{\hy{#1}{\Al}}
\newcommand{\Abgp}{\mbox{\sf Abgp}}
\newcommand{\Ap}{{\EuScript A}\sb{p}}
\newcommand{\Arr}{\mbox{\sf Arr}}
\newcommand{\Cat}{\mbox{\sf Cat}}
\newcommand{\DGLA}{\mbox{\sf DGLA}}
\newcommand{\Map}{\mbox{\sf Map}}
\newcommand{\sMap}{\mbox{\sf sMap}}
\newcommand{\dsMst}{\sMap\sp{\Stov,R}}
\newcommand{\dsMsr}{\sMap\sb{\real}\sp{\Stov,R}}
\newcommand{\sMst}{\sMap\sb{\Stov}}
\newcommand{\sMstn}[1]{\sMap\sb{\Stov}\sp{#1}}
\newcommand{\wMap}{\mbox{\sf wMap}}
\newcommand{\Set}{\mbox{\sf Set}}
\newcommand{\Seta}{\Set\sb{\ast}}
\newcommand{\nSys}[1]{\hy{#1}{\Sys}}
\newcommand{\Grp}{\mbox{\sf Gp}}
\newcommand{\Top}{\mbox{\sf Top}}
\newcommand{\Ta}{\Top\sb{\ast}}
\newcommand{\Tz}{\Top\sb{0}}
\newcommand{\RM}[1]{{#1}\text{-}{\EuScript Mod}}
\newcommand{\gn}[1]{\gamma\sb{[#1]}}
\newcommand{\Sa}{\Ss\sb{\ast}}
\newcommand{\Sn}[1]{\Ss\sb{[#1]}}
\newcommand{\Sk}{\Ss\sp{\Kan}}
\newcommand{\Sr}{\Ss\sp{\red}}
\newcommand{\SR}{\Ss\sb{R}}
%
%
\newcommand{\FF}{\mathbb F}
\newcommand{\Fp}{\FF\sb{p}}
\newcommand{\NN}{\mathbb N}
\newcommand{\QQ}{\mathbb Q}
\newcommand{\ZZ}{\mathbb Z}
%
%
\newcommand{\bA}{{\mathbf A}}
\newcommand{\bDel}{{\mathbf \Delta}}
\newcommand{\bDeu}{\bDel\sp{\bullet}}
\newcommand{\Dz}{\Delta\sb{+}}
\newcommand{\Dop}{\Delta\op}
\newcommand{\Dzop}{\Dz\op}
\newcommand{\bpar}[1]{{\mathbf\partial}\sb{#1}}
\newcommand{\be}[1]{{\mathbf e}\sp{#1}}
\newcommand{\gS}[1]{{\EuScript S}\sp{#1}}
%
%
\newcommand{\Fn}{\{K(\Fp,i)\}\sb{i=1}\sp{\infty}}
\newcommand{\KR}[1]{K(R,{#1})}
\newcommand{\Rn}{\{\KR{n}\}\sb{n=1}\sp{\infty}}
%
%
\newcommand{\bB}{{\mathbf  B}}
\newcommand{\bD}{{\mathbf  D}}
\newcommand{\sD}{{\mathbf \Sigma D}}
\newcommand{\bLL}{{\mathbf L}}
\newcommand{\bLr}{\bLL\sp{\rel}}
\newcommand{\bM}{{\mathbf M}}
\newcommand{\bo}{\mathbf 1}
\newcommand{\bRR}{{\mathbf R}}
\newcommand{\bRr}{\bRR\sp{\rel}}
\newcommand{\bS}[1]{{\mathbf S}\sp{#1}}
\newcommand{\bU}{{\mathbf U}}
\newcommand{\bW}{{\mathbf W}}
\newcommand{\bX}{{\mathbf X}}
\newcommand{\bY}{{\mathbf Y}}
\newcommand{\bZ}{{\mathbf Z}}
%
%
\newcommand{\pis}{\pi\sb{\ast}}
\newcommand{\pin}[1]{\pi\sp{\natural}\sb{#1}}
\newcommand{\piA}{\pis\sp{\A}}
\newcommand{\piS}{\pis\sp{\uS}}
%
%
\newcommand{\Eut}[2][ ]{$E\sp{#2}$-term{#1}}
\newcommand{\Elt}[2][ ]{$E\sb{#2}$-term{#1}}
%
%
\newcommand{\Pa}[1][ ]{$\Pi$-algebra{#1}}
\newcommand{\PiS}{\Pi\sb{\uS}}
\newcommand{\PSa}[1][ ]{$\PiS$-algebra{#1}}
\newcommand{\PAlg}{\Alg{\Pi}}
\newcommand{\PSAlg}{\Alg{\PiS}}
\newcommand{\PiA}{\Pi\sp{\A}}
\newcommand{\PAa}[1][ ]{$\PiA$-algebra{#1}}
\newcommand{\PAAlg}{\Alg{\PiA}}
%
%
\newcommand{\co}[1]{c({#1})\sb{\bullet}}
\newcommand{\cu}[1]{c({#1})\sp{\bullet}}
\newcommand{\Ad}{A\sb{\bullet}}
\newcommand{\Au}{A\sp{\bullet}}
\newcommand{\As}{A\sb{\ast}}
\newcommand{\Bd}{B\sb{\bullet}}
\newcommand{\Bu}{B\sp{\bullet}}
\newcommand{\Gd}{G\sb{\bullet}}
\newcommand{\oG}[1]{\overline{G}\sb{#1}}
\newcommand{\Qd}{Q\sb{\bullet}}
\newcommand{\Ud}{U\sb{\bullet}}
\newcommand{\Uu}{\bU\sp{\bullet}}
\newcommand{\Vd}{V\sb{\bullet}}
\newcommand{\oV}[1]{\overline{V}\sb{#1}}
\newcommand{\Vu}{V\sp{\bullet}}
\newcommand{\Wd}{\bW\sb{\bullet}}
\newcommand{\Yd}{Y\sb{\bullet}}
\newcommand{\Zd}{\bZ\sb{\bullet}}
%
%
\newcommand{\irn}[1]{\iota\bp{#1}}
\newcommand{\irnk}[2]{\irn{#1}\sb{#2}}
\newcommand{\bv}{{\bm \vare}}
\newcommand{\bet}{{\bm \eta}}
\newcommand{\bve}[1]{\bv\bp{#1}}
\newcommand{\cW}{{\EuScript W}}
\newcommand{\Wn}[2]{\bW\sb{#1}\bp{#2}}
\newcommand{\WW}[1]{\Wn{\bullet}{#1}}
\newcommand{\tW}{\widetilde{\bW}}
\newcommand{\tWn}[2]{\tW\quad\hspace*{-4mm}\sb{#1}\bp{#2}}
\newcommand{\tWd}[1]{\tWn{\bullet}{#1}}
\newcommand{\btW}{{\raisebox{-1.5ex}{$\stackrel{\textstyle\overline{\bW}}{\sim}$}}}
\newcommand{\btWn}[2]{\btW\quad\hspace*{-4mm}\sb{#1}\bp{#2}}
\newcommand{\oW}[1]{\overline{\bW}\sb{#1}}
\newcommand{\ppp}{\hspace*{0.3mm}\sp{\prime}\hspace{-0.3mm}}
\newcommand{\oWp}[1]{\ppp\oW{#1}}
\newcommand{\Wpn}[2]{\ppp\bW\sb{#1}\bp{#2}}
\newcommand{\Wp}[1]{\Wpn{\bullet}{#1}}
\newcommand{\oWn}[2]{\oW{#1}\bp{#2}}
\newcommand{\uW}[1]{\overline{\bW}\hspace{0.2mm}\sp{#1}}
\newcommand{\hW}[1]{\widehat{\bW}\sb{#1}}
\newcommand{\hWp}[1]{\ppp\hW{#1}}
\newcommand{\vW}[1]{\widehat{\bW}\sb{#1}}
\newcommand{\vWn}[2]{\vW{#1}\bp{#2}}
\newcommand{\vWu}[1]{\vWn{\bullet}{#1}}
\newcommand{\dvW}[1]{\widehat{\bW}\sp{#1}}
\newcommand{\dvWn}[2]{\dvW{#1}\lp{#2}}
\newcommand{\dvWu}[1]{\dvWn{\bullet}{#1}}
\newcommand{\SW}[2]{\Sigma\sp{#1}\oW{#2}}
\newcommand{\CsW}[2]{C\Sigma\sp{#1}\oW{#2}}
\newcommand{\Xd}{\bX\sb{\bullet}}
\newcommand{\cZ}{{\EuScript Z}}
%
%
\newcommand{\Du}[1]{\bD\sp{\bullet}\lp{#1}}
\newcommand{\sDu}[1]{\sD\sp{\bullet}\lp{#1}}
\newcommand{\Gu}{G\sp{\bullet}}
\newcommand{\prn}[1]{p\lp{#1}}
\newcommand{\prnk}[2]{\prn{#1}\sp{#2}}
\newcommand{\Wu}{\bW\sp{\bullet}}
\newcommand{\Wun}[2]{\bW\sp{#1}\lp{#2}}
\newcommand{\WWu}[1]{\Wun{\bullet}{#1}}
\newcommand{\ooW}[2]{\Omega\sp{#1}\uW{#2}}
\newcommand{\PoW}[2]{P\Omega\sp{#1}\uW{#2}}
\newcommand{\tWun}[2]{\tW\quad\hspace*{-4mm}\sp{#1}\lp{#2}}
\newcommand{\tWu}[1]{\tWun{\bullet}{#1}}
\newcommand{\vWun}[1]{\vW{\bullet}\bp{#1}}
\newcommand{\Xu}{\bX\sp{\bullet}}
%
%
\newcommand{\Ama}[1][ ]{dual $\uA$-mapping algebra{#1}}
\newcommand{\sAma}[1][ ]{dual strict $\uA$-mapping algebra{#1}}
\newcommand{\dsStma}[1][ ]{dual strict Stover mapping algebra{#1}}
\newcommand{\wAma}[1][ ]{dual weak $\uA$-mapping algebra{#1}}
\newcommand{\eAma}[1][ ]{extended $\A$-mapping algebra{#1}}
\newcommand{\eApma}[1][ ]{extended $\A'$-mapping algebra{#1}}
\newcommand{\hAma}[1][ ]{$\hA$-mapping algebra{#1}}
\newcommand{\Fma}[1][ ]{$\Fp$-mapping algebra{#1}}
\newcommand{\sRma}[1][ ]{dual strict $R$-mapping algebra{#1}}
\newcommand{\wRma}[1][ ]{dual weak $R$-mapping algebra{#1}}
\newcommand{\sMa}[1][ ]{strict mapping algebra{#1}}
\newcommand{\Sma}[1][ ]{$\uS$-mapping algebra{#1}}
\newcommand{\sSma}[1][ ]{strict $\uS$-mapping algebra{#1}}
\newcommand{\Stma}[1][ ]{Stover mapping algebra{#1}}
\newcommand{\sStma}[1][ ]{strict Stover mapping algebra{#1}}
\newcommand{\wSma}[1][ ]{weak $\uS$-mapping algebra{#1}}
\newcommand{\bT}{\mathbf{T}}
\newcommand{\bTh}{\mathbf{\Theta}}
\newcommand{\TsS}{\bTh\sb{\uS}}
\newcommand{\ThsS}{\Theta\sb{\uS}}
\newcommand{\TsSt}{\bTh\sb{\Stov}}
\newcommand{\TuA}{\bTh\sp{\uA}}
\newcommand{\ThuA}{\Theta\sp{\uA}}
\newcommand{\TuR}{\bTh\sp{R}}
\newcommand{\TRl}{\TuR\sb{\lambda}}
\newcommand{\TuSR}{\bTh\sp{\Stov,R}}
\newcommand{\TuSRl}{\TuSR\sb{\lambda}}
\newcommand{\ma}[1][ ]{mapping algebra{#1}}
\newcommand{\Rma}[1][ ]{$R$-mapping algebra{#1}}
%
%
\newcommand{\MuA}{\Map\sp{A}}
\newcommand{\sMuA}{ \sMap\sp{A}}
\newcommand{\wMuA}{ \wMap\sp{A}}
\newcommand{\pMuA}{'\hspace{-0.5mm}\Map\sp{A}}
\newcommand{\sMuR}{ \sMap\sp{R}}
\newcommand{\MsS}{\Map\sb{\uS}}
\newcommand{\sMsS}{\sMap\sb{\uS}}
\newcommand{\sMsSn}[1]{\sMap\sb{\uS}\sp{#1}}
\newcommand{\sMsSt}{\sMap\sb{\Stov}}
\newcommand{\wMsS}{\wMap\sb{\uS}}
\newcommand{\wMsSn}[1]{\wMap\sb{\uS}\sp{#1}}
\newcommand{\fff}{\mathfrak{f}}
\newcommand{\tff}[1]{\widetilde{\fff}\bp{#1}}
\newcommand{\hff}[1]{\widehat{\fff}\bp{#1}}
\newcommand{\tf}[2]{\tilde{f}\bp{#2}\sb{#1}}
\newcommand{\of}[2]{\bar{f}\up{#1}\sb{#2}}
\newcommand{\fg}{\mathfrak{g}}
\newcommand{\fM}{\mathfrak{M}}
\newcommand{\hfM}{\hat{\fM}}
\newcommand{\fMA}{\fM\sp{\uA}}
\newcommand{\fMR}{\fM\sp{R}}
\newcommand{\fMS}{\fM\sb{\uS}}
\newcommand{\fMst}{\fM\sb{\Stov}}
\newcommand{\fMRst}{\fM\sp{\Stov,R}}
\newcommand{\fT}{\mathfrak{T}}
\newcommand{\fTn}[1]{\fT\sp{#1}}
\newcommand{\fU}{\mathfrak{U}}
\newcommand{\fUd}{\fU\sb{\bullet}}
\newcommand{\fV}{\mathfrak{V}}
\newcommand{\fVd}{\fV\sb{\bullet}}
\newcommand{\fW}{\mathfrak{W}}
\newcommand{\fWd}{\fW\sb{\bullet}}
\newcommand{\fX}{\mathfrak{X}}
\newcommand{\hfX}{\hat{\fX}}
\newcommand{\fXd}{\fX\sb{\bullet}}
\newcommand{\fY}{\mathfrak{Y}}
\newcommand{\fYd}{\fY\sb{\bullet}}
\newcommand{\fZ}{\mathfrak{Z}}
\newcommand{\bone}{[\mathbf{1}]}
\newcommand{\bbn}{[\mathbf{n}]}
\newcommand{\bbk}{[\mathbf{k}]}
\setcounter{section}{-1}
%
%
\title{Truncated Derived Functors and Spectral Sequences}
\date{\today}
\author{Hans-Joachim Baues}
\address{Max-Planck-Institut f\"{u}r Mathematik\\
Vivatsgasse 7\\ 53111 Bonn, Germany}
\email{baues@mpim-bonn.mpg.de}
\author{David Blanc}
\address{Department of Mathematics\\ University of Haifa\\ 3498838 Haifa\\ Israel}
\email{blanc@math.haifa.ac.il}
\author{Boris Chorny}
\address{Department of Mathematics\\ University of Haifa at Oranim\\ 3600600 Tivon\\ Israel}
\email{chorny@math.haifa.ac.il}
\date{\today}

\subjclass[2010]{Primary: 55T99; \ secondary: 55U35, 18G50, 18C30}

\keywords{Spectral sequence, relative derived functors, truncation, differentials,
(co)simplicial resolutions, mapping algebra}

\begin{abstract}
  The \Elt{2} of the Adams spectral sequence may be identified with certain derived
  functors, and this also holds for a number of other spectral sequences.
  Our goal is to show how the higher terms of such spectral sequences are determined
  by truncations of relative derived functors, defined in terms of certain simplicial
functors called \emph{mapping algebras}.
\end{abstract}

\maketitle
%
%
\sect{Introduction}

The various types of Adams spectral sequences, which play a central role in
algebraic topology (cf.\ \cite{AdSS,BCMiU,BCurS,BKanS,BKanH,NovO,RavC},
have a number of features in common:

\begin{enumerate}
\renewcommand{\labelenumi}{(\roman{enumi})}
\item They are obtained from a space $\bY$ by
  constructing a (cosimplicial) resolution \w{\bY\to\Wu} with respect to a spectrum
  \w[,]{\uA=\{A\sb{i}\}\sb{i=-\infty}\sp{\infty}} with its associated
  cohomology theory \w[.]{\uA\sp{\ast}}
\item The spectral sequence in question is the homotopy spectral sequence for
\w[,]{\bT\Wu} for a suitable homotopy functor $\bT$.
\item The \Elt{2} of the spectral sequence can be identified as the derived functors
of an algebraic functor $T$ associated to $\bT$, applied to \w[.]{\uA\sp{\ast}\bY}
\end{enumerate}

The goal of this paper is to provide a description similar to (iii) for the \Elt{n+2}
of the spectral sequence (for \w[),]{n\geq 0} as relative derived functors
applied to the truncation \w{\Po{n}\fMA\bY} of a certain structure, called a
\emph{\ma[,]} associated to $\bY$ (which reduces to \w{\uA\sp{\ast}\bY} when \w[).]{n=0}

Just as for the \Elt[,]{2} this has two advantages:
\begin{enumerate}
\renewcommand{\labelenumi}{(\alph{enumi})}
\item The truncated \ma \w{\Po{n}\fMA\bY} has less information than $\bY$ itself,
but still enough to determine the \Elt[.]{n+2}
\item Relative derived functors may be calculated using \emph{any} resolution
  of \w[.]{\Po{n}\fMA\bY}
\end{enumerate}

The first author carried out this program for the \Elt{3} of the stable Adams
spectral sequence in \cite{BauAS,BJiblD}, showing that extended calculations may be
made using such a construction. See \cite{BBlaH,CFranH} for other general descriptions
of the higher terms in the stable Adams spectral sequence, although not quite in
the form of truncated derived functors as defined here.

\begin{mysubsection}{Mapping algebras and truncations}
\label{smapalg}
By (iii) above, the \Elt{2} of the Adams spectral sequence
depends only on the sets \w{[\bY,A\sb{i}]\sb{i\in\ZZ}} and
operations on them induced by homotopy classes of maps between (products of) the
spaces \w[.]{A\sb{i}} This suggests that for the higher terms, we should look at
the function spaces \w[,]{\mapa(\bY,A\sb{i})} with additional structure induced by maps
between the representing spaces. This structure is encoded by the notion of
a \emph{\ma[:]} that is, a simplicial functor \w{\fX:\TuA\to\Sa} from the
sub-simplicial category \w{\TuA} of \w{\Tz} whose objects are
products of copies of the various spaces \w[.]{A\sb{i}}
For example, the \emph{realizable} \ma \w{\fX:=\fMA\bY} has the
value \w{\map(\bY,\bA)} at each \w[.]{\bA\in\TuA}

Mapping algebras admit truncations, defined by applying the Postnikov
section functor \w{\Po{n}} to each mapping space. In particular, the $0$-truncation
contains the same information as the sets \w{[\bY,A\sb{i}]\sb{i\in\ZZ}} of homotopy
classes of maps, together with the operations on them induced by homotopy classes of
maps between the spaces \w[:]{A\sb{i}} this is precisely what was needed to determine
the \Eut{2} as suitable derived functors in (iii) above.

This suggests that higher truncations of the mapping algebras may suffice
to determine higher terms in the spectral sequence \wh depending, of course, on
the homotopy functor $\bT$ in question.

We may therefore summarize our program as follows:

\begin{enumerate}
\renewcommand{\labelenumi}{(\arabic{enumi})}
\item We need to show how a continuous functor \w{\bT:\Ta\to\Ta}
  factors through the category \w{\MuA} of \ma[s] as
  \w[,]{\fT\circ\fMA} for a suitable homotopy functor \w[.]{\fT:\MuA\to\Ta}
\item We want \w{\fWd:=\fMA\Wu} to be a resolution of
  \w{\fMA\bY} in the resolution model category of simplicial \ma[s,] in order to
  guarantee that both the (functorial) cosimplicial resolution \w{\Wu} of $\bY$,
  and the resulting cosimplicial space \w[,]{\bT\Wu} are homotopy functors of
  \w[.]{\fMA\bY} This will let us identify \w{\bT\Wu} as a certain relative
  left derived functor \w{(\bLr\fT)\fMA\bY=\fT\fWd} of $\fT$
  applied to the \ma \w{\fMA\bY} (see \S \ref{srdf}).
\item Finally, we must show that in the cases of interest to us, the \Elt{r} of the
  homotopy spectral sequence for \w{\bT\Wu=(\bLL\fT)\fMA\bY} depends only on the
  $n$-truncation \w[,]{\Po{r+2}\fWd} for each \w[.]{r\geq 2} Functors $\bT$ with this
  property are called \emph{level}.
\end{enumerate}
\end{mysubsection}

\begin{remark}\label{rdualvers}
There are also a number of less familiar spectral sequences obtained dually
by constructing a simplicial resolution \w{\Xd\to\bY} with respect to
\w[,]{\uS=\{\bS{i}\}\sb{i=1}\sp{\infty}} applying a homotopy functor
\w[,]{\bT:\C\to\C} and then using the homotopy spectral sequence for
the simplicial space\w{\bT\Xd} (see \cite{StoV,BlaH,DKSStP}).  Here too, one can identify
the \Eut{2} with the derived functors of an algebraic functor of \w{\pis\bY}
(the algebraic object corepresented by $\uS$). We include these in the paper
mainly in order to show that the formalism we describe here is not limited to
the Adams spectral sequence, even though this is our most important example.
Moreover, in a number of ways the simplicial-covariant version is cleaner than
the cosimplicial-contravariant one.

However, since Eckmann-Hilton duality is not formal, we are forced to work carefully
through the details in the two versions separately: for this reason, each section is
divided into two parts, starting with the covariant case.

For reasons of space, we deal here only with the unstable spectral sequences.
For the stable analogue, we must choose a simplicial model category of spectra
(cf.\ \cite{BFrieH,EKMMayR,HSSmiS,LydaS})
and work there throughout; one can still take Postnikov $n$-sections of the
mapping spaces \w[.]{\mapa(\uS,\,\Xd)}
\end{remark}

\begin{mysubsection}{Outline}
\label{sorg}
  In Section \ref{cmapalg} we define enriched sketches and the associated
  \ma[s] (as well as the dual versions). It turns out that we have competing
  versions of \ma[s:] the category \w[,]{\dsMsr} which allows us to
  factor $\bT$ as \w{\fT\circ\fMA} in \S \ref{smapalg}(1), is not right proper,
  so we need a variant \w{\Sa\sp{\TuA}} in which \w{\fWd:=\fMA\Wu} is
indeed a cofibrant replacement for \w{\fMA\bY} in the  resolution model category
  \w[.]{\Sa\sp{\TuA\times\Dop}}

  In Section \ref{cffma} we construct the category \w{\dsMsr} of \ma[s,] for
  $\cA$ the Eilenberg-Mac~Lane spectrum for a commutative ring $R$, and show:
\begin{thma}
  There is a \emph{realization} functor
  \w[,]{N:(\dsMsr)\op\to\Sa} equipped with a natural weak equivalence
  \w[.]{\fMA\circ N\to\Id}
\end{thma}
\noindent See Theorem \ref{tdrssma} and Corollary \ref{cdrssma} below.

Thus any homotopy functor \w{\bT:\Ta\to\Ta} which preserves $R$-equivalences,
when restricted to $R$-good spaces, induces a functor \w{\fT:=\bT\circ N:(\dsMsr)\op\to\Ta}
equipped with a natural weak equivalence \w[.]{\fT\circ\fMRst\to\bT}

When $\uS$ is the sphere spectrum (cf.\ \S \ref{rdualvers}), there is a dual category
\w{\sMst} of \Sma[s] with a realization functor \w[,]{N:\sMst\to\Ta}
(see Theorem \ref{trssma} and Corollary \ref{crssma}).

In Section \ref{chotdf} we define the general notion of a relative derived functor
(\S \ref{srdf}), and show how it applies to the functor
\w{\fT:(\dsMsr)\op\to\Ta} associated to the homotopy functor \w[.]{\bT:\Ta\to\Ta}
To do so, we have to relate the two types of \ma[s] described in
Section \ref{cmapalg} \wh those that are used for resolutions,
and those for which $\fT$ is defined \wh by means of Theorem \ref{tdtlwe}, which implies:
\begin{thmb}
  If $\bY$ is $R$-good, any simplicial resolution \w{\fVd} of \w{\fMRst\bY} in the
  resolution model category \w{\Sa\sp{\TuA\times\Dop}} is Reedy weakly equivalent
  (i.e., in each simplicial dimension) to a simplicial object \w{\fWd} in
  \w[.]{(\dsMsr)\sp{\Dop}}
\end{thmb}
\noindent The dual version, for the sphere spectrum, is Theorem \ref{ttlwe}.

Finally, in Section \ref{cthodf} we deal with the truncated versions of our higher
derived functors, explain what data is needed to determine the \Elt{r} of the homotopy
spectral sequence of a (co)simplicial space by formalizing the notion of a \emph{level
functor} (\S \ref{dlevelf}), and show

\begin{thmc}
For \w{R=\Fp} or $\QQ$, \w[,]{\bZ\in\Sa} and $R$-good $\bY$, the unstable
Adams spectral sequence for \w{\mapa(\bZ,\bY)} is determined by a
simplicial \ma resolution \w{\fWd} of \w[,]{\fMRst\bY} and for each \w{r\geq 2}
the \Elt{r} is determined by the corresponding \wwb{r-2}truncated \ma[s.]
\end{thmc}
\noindent See Theorem \ref{tdlevelf}.

This implies that the mapping space functor \w{\mapa(\bZ,-)} is a level homotopy
functor on $R$-good spaces.
We also prove a number of similar results for functors related to the sphere spectrum
(see Propositions \ref{phure}, \ref{psusp}, and \ref{pwedge}).
\end{mysubsection}

\begin{notn}\label{snac}
The category of finite ordered sets and order-preserving maps will be denoted
by $\Delta$ (cf.\ \cite[\S 2]{MayS}), so a simplicial object
\w{\Gd} in $\C$ is a functor \w[,]{\Dop\to\C} and the category of such will
be denoted by \w[.]{\C\sp{\Dop}} Similarly, a cosimplicial object \w{\Gu}
in a category $\C$ is a functor \w[,]{\Delta\to\C} and the category of such will
be denoted by \w[.]{\C\sp{\Delta}}
There is a natural embedding \w{\co{-}:\C\to\C\sp{\Dop}} (the constant simplicial
object), and similarly \w[.]{\cu{-}:\C\to\C\sp{\Delta}}

Write \w{\Delta\sb{+}} for the subcategory
of $\Delta$ with the same objects but only  monic maps. A functor \w{G:\Dzop\to\C}
(respectively, \w[)]{G:\Dz\to\C} is called a \emph{restricted} (co)simplicial
object in $\C$. The inclusion \w{i:\Dz\hra\Delta} induces a forgetful functor
\w[,]{i\sp{\ast}:\C\sp{\Dop}\to \C\sp{\Dzop}} which has a left adjoint
\w{\LL:\C\sp{\Dzop}\to \C\sp{\Dop}} (for suitable $\C$).

The category of topological spaces will be denoted by \w[,]{\Top}
that of pointed spaces by \w[,]{\Ta} and that of pointed connected
spaces by \w[.]{\Tz} The category of simplicial sets will be
denoted by \w[,]{\Ss=\Set\sp{\Dop}} that of pointed simplicial
sets by \w[,]{\Sa=\Seta\sp{\Dop}} that of simplicial groups by
\w[.]{\G=\Grp\sp{\Dop}} Write \w{\mapa(\bX,\bY)} for the standard function complex in
\w[,]{\Sa} \w[,]{\Tz} or $\G$ (see \cite[I, \S 1.5]{GJarS}).
Note that both \w{\Tz} and \w{\Sa} are enriched over \w[,]{(\Sa,\wedge)}
but if we forget the basepoints, the same mapping spaces \w{\map\sb{\Sa}(X,Y)} or
\w{\map\sb{\Tz}(X,Y)} also define an enrichment over \w[,]{(\Ss,\times)} which is
the one we shall use (see \cite[9.1.14]{PHirM}),

We denote the category of pointed Kan complexes by
\w[,]{\Sk} that of \emph{reduced} simplicial sets (with a single
vertex) by \w[,]{\Sr} and the full subcategory of $n$-types in
\w{\Sa} \wh i.e., spaces $X$ with \w{\pi\sb{i}(X,x)=0} for \w{i>n}
and all \w{x\in X\sb{0}} \wwh by \w[,]{\Sn{n}}  with \w{\Po{n}:\Sa\to\Sn{n}}
the $n$-th \emph{Postnikov section} functor.
\end{notn}

\begin{ack}
The research of the second and third authors was partially
supported by Israel Science Foundation grants 770/16 and 1138/16,
respectively.
\end{ack}

%
%
\sect{Mapping algebras}
\label{cmapalg}

The main technical tool in our approach is the notion of a \ma[,] first
introduced in \cite[\S 8]{BBlaC}.  We shall need a number of variants of this notion,
together with their dual versions\vsm.

\supsect{\protect{\ref{cmapalg}}.A}{Enriched sketches and mapping algebras}

\begin{defn}\label{dmapalg}
Let $\C$ be a pointed simplicial model category, $\uS$ a set of fibrant and cofibrant
homotopy cogroup objects in $\C$, $\F$ a category of finite simplicial sets, and $\E$ a
set of cocones in $\C$. The associated \emph{enriched sketch}, or multi-sorted theory
(cf.\ \cite[\S 5.6]{BorcH2})
\w{\TsS=\bTh\sb{(\uS,\F,\E)}} is the smallest full sub-simplicial category of
$\C$ containing $\uS$ and closed under the operations \w{-\otimes K} for
\w{K\in\F} and taking colimits of the cocones in $\E$. In this setting:
\begin{enumerate}
\renewcommand{\labelenumi}{(\arabic{enumi})}
\item A \emph{$\uS$-presheaf} is a pointed simplicial functor
  \w[.]{\fX:\TsS\op\to\Sa}  The category of all $\uS$-presheaves is denoted by
  \w[,]{\Sa\sp{\TsS\op}} and the value of $\fX$ at \w{\bB\in\TsS} will be written
  \w[.]{\fX\lin{\bB}}

A map \w{\fff:\fX\to\fY} of $\uS$-presheaves is called a \emph{weak equivalence} if
\w{\fff\lin{\bB}:\fX\lin{\bB}\to\fY\lin{\bB}} is a weak equivalence for each
\w[.]{\bB\in\TsS}
Two $\uS$-presheaves are said to be \emph{weakly equivalent} if they are connected by
a finite zigzag of weak equivalences.
\item A \emph{\sSma} is a $\uS$-presheaf $\fX$ for which the natural maps
\begin{myeq}\label{eqsuscoprod}
\fX\lin{\bB\otimes K}~\to~\fX\lin{\bB}\sp{K}~\hspace{9mm}\text{and}\hspace{9mm}
\fX\lin{\colim\sb{i\in I}\,\bB\sb{i}}~\to~\lim\sb{i\in I}\,\fX\lin{\bB\sb{i}}
\end{myeq}
\noindent are isomorphisms for all \w[,]{\bB\in\TsS} \w[,]{K\in\F} and
diagrams $I$ in $\E$. The full subcategory of \sSma[s] will be denoted by \w[.]{\sMsS}
\item A \emph{\wSma} is a $\uS$-presheaf $\fX$ which is weakly equivalent to a \sSma[.]
Thus in particular, the maps of \wref{eqsuscoprod} are weak equivalences.
The full subcategory of \wSma[s]  will be denoted by \w[.]{\wMsS}
\end{enumerate}
\end{defn}

\begin{remark}\label{rwsma}
In principle, we would like to identify a \wSma more conceptually as a
$\uS$-presheaf for which not only the maps of \wref{eqsuscoprod} are weak equivalences,
but also appropriate higher coherences hold. However, as we shall not in fact need
to work explicitly with \wSma[s,] we can make do here with the above ad hoc definition.
\end{remark}

\begin{example}\label{egmapalg}
The main example of an enriched sketch we shall consider in this paper is
the case where \w[,]{\C=\Tz} \w{\uS=\{\bS{n}\}\sb{n=1}\sp{\infty}} and $\F$
consists of the inclusions \w[,]{i\sb{0},i\sb{1}:\Delta[0]\hra\Delta[1]}
The cocone collection $\E$ contains all coproducts of cardinality $<\lambda$
for some fixed limit cardinal $\lambda$ (e.g., \w[),]{\aleph\sb{0}} and the
pushout squares
\mydiagram[\label{eqconesusp}]{
\ar @{} [drr] |(0.7){\framebox{\scriptsize{PO}}}
\bB \ar@{^{(}->}[rr] \ar@{->>}[d]\sb{\simeq} &&
\bB\otimes\Delta[1] \ar@{->>}[d]\sb{\simeq}\sp{\inc\sb{0}}&&
\ar @{} [drr] |(0.7){\framebox{\scriptsize{PO}}}
\bB \ar@{^{(}->}[rr] \ar@{->>}[d] && C\bB \ar@{->>}[d]\sp{\inc\sb{1}} \\
\ast \ar@{^{(}->}[rr] && C\bB && \ast \ar@{^{(}->}[rr] && \Sigma\bA
}
\noindent for \w[.]{\bB\in\TsS} (These will be our models for the cone \w{CX} and
suspension \w{\Sigma X} of any \w[).]{X\in\C}

Thus a \sSma $\fX$ will take the two squares of \wref{eqconesusp}
to pullback squares:
\mydiagram[\label{eqcotolim}]{
\ar @{} [drr] |(0.3){\framebox{\scriptsize{PB}}}
P\fX\lin{\bB} \ar@{^{(}->}[rr] \ar@{->>}[d]\sb{\simeq} &&
\fX\lin{\bB}\sp{\Delta[1]} \ar@{->>}[d]\sp{\simeq}\sb{\ev\sb{0}}&&
\ar @{} [drr] |(0.3){\framebox{\scriptsize{PB}}}
\Omega\fX\lin{\bB} \ar@{^{(}->}[rr]\sp{\iota\sb{\bB}} \ar@{->>}[d] &&
P\fX\lin{\bB} \ar@{->>}[d]\sb{\ev\sb{1}} \\
\ast \ar@{^{(}->}[rr] && \fX\lin{\bB} && \ast \ar@{^{(}->}[rr] && \fX\lin{\bB}
}

One might also consider localized versions, where
\w{\uS=\{\bS{n}\sb{R}\}\sb{n=1}\sp{\infty}} for some subring \w{R\subseteq\QQ}
(cf.\ \cite{BittS}).  In particular, when \w{R=\QQ} we
may replace \w{\C=\Tz} by a suitable algebraic model of rational homotopy types, such as
the category of differential graded Lie algebras.

More generally, one could take any space \w[,]{\bM\in\Tz} and let
\w[.]{\uS=\{\Sigma\sp{n}\bM\}\sb{n=1}\sp{\infty}} However, while the formal part of our
program can be made to work in this case (see \cite{BBDoraR}), the application to
the homotopy spectral sequence of a simplicial space is not available for $\bM$ which
is not essentially a sphere (see \cite{CDInteA} and \cite[\S 4.6]{BlaM}).
\end{example}

\begin{defn}\label{drdma}
For any enriched sketch \w{\TsS} as above, the most important example of a
$\uS$-presheaf $\fX$ is a \emph{realizable} one, associated to an object
\w[,]{\bY\in\C}  where \w{\fX\lin{\bB}:=\map\sb{\C}(\bB,\bY)} for any \w[.]{\bB\in\TsS}
Evidently, this will be a \sSma[,] which we denote by \w{\fMS\bY} (of course,
it actually takes \emph{all} colimits in \w{\TsS} to the corresponding limits).
When \w[,]{\bY\in\Obj\TsS} we say that \w{\fMS\bY} is \emph{free}.
\end{defn}

The strong Yoneda Lemma for enriched categories (see \cite[2.4]{GKellyEC}) implies:

\begin{lemma}\label{lfreema}
If $\fY$ is a $\uS$-presheaf and \w{\fMS\bB} is a free \sSma (for \w[),]{\bB\in\TsS}
there is a natural isomorphism
$$
\Phi:\map\sb{\Sa\sp{\TsS\op}}(\fMS\bB,\,\fY)~\xra{\cong}~\fY\lin{\bB}~,
$$
\noindent with \w{\Phi(\fff)=\fff(\Id\sb{\bB})\in\fY\lin{\bB}\sb{0}} for any
\w[.]{\fff\in\Hom\sb{\Sa\sp{\TsS\op}}(\fMS\bB,\,\fY)=
  \map\sb{\Sa\sp{\TsS\op}}(\fMS\bB,\,\fY)\sb{0}}
\end{lemma}

\begin{remark}\label{raltma}
It is sometimes convenient think of a $\uS$-presheaf $\fX$ as
a category $\cX$ with object set
\w[,]{\OO:=\Obj(\TsS)\cup\{\star\}} enriched in pointed simplicial sets as follows:
\begin{myeq}\label{eqaltma}
\map\sb{\cX}(\bA,\bB)~=~\begin{cases}
\map\sb{\TsS}(\bA,\bB) & \hsm\text{if}\hsm \bA,\bB\in\Obj(\TsS)\\
\fX\lin{\bA} & \hsm\text{if}\hsm \bA\in\Obj(\TsS)\hsm \text{and}\hsm \bB=\star\\
\co{\{\ast,\Id\sb{\bstar}\}} & \hsm\text{if}\hsm \bA=\bB=\bstar\\
\co{\{\ast\}} & \hsm\text{otherwise.}
\end{cases}
\end{myeq}
\noindent Thus a realizable $\uS$-presheaf \w{\fX=\fMS\bY} corresponds to a sub-simplicial
category $\cX$ of $\C$ with object set \w{\Obj(\TsS)\cup\{\bY\}}
(compare \cite[\S 8.1]{BBlaC}).
\end{remark}

\begin{defn}\label{dpalg}
An enriched sketch \w{\TsS} in a model category $\C$ has an algebraic version,
which is the (ordinary) sketch \w{\ThsS:=\pi\sb{0}\TsS} \wwh
that is, \w{\ThsS} has the same objects as \w[,]{\TsS} and
\w[.]{\Hom\sb{\ThsS}(\bB,\bB'):=\pi\sb{0}\map\sb{\TsS}(\bB,\bB')}
An \emph{algebra} (or \emph{model}) for \w{\ThsS} is a functor \w{\Lambda:\ThsS\op\to\Set}
which takes the coproduct of any discrete cocone in $\E$ to a product in \w{\Set}
(see \cite[\S 5.6]{BorcH2}).

These are called \emph{\PSa[s]}, and the category of such is denoted
by \w[:]{\PSAlg} for \w[,]{\uS=\{\bS{n}\}\sb{n=1}\sp{\infty}} these are simply
the \emph{\Pa[s]} of \cite{DKanE}.
Note that if $\fX$ is a (weak or strict) \Sma[,] then \w{\pi\sb{0}\fX} is a \PSa[;]
the same need not hold for an arbitrary $\uS$-presheaf.  We say that a \PSa $\Lambda$
is \emph{realizable} if it is of the form \w{\pi\sb{0}\fMS\bY} for some
\w[.]{\bY\in\C} A coproduct of \PSa[s] of the form
\w{\pi\sb{0}\fMS\bB} for \w{\bB\in\Obj\TsS} is called \emph{free}.
\end{defn}

\supsect{\protect{\ref{cmapalg}}.B}{Dual sketches and mapping algebras}

There are dual versions of all three notions, defined as follows:

\begin{defn}\label{dsmapalg}
  Let $\C$ be a pointed simplicial model category, $\uA$ a set of fibrant and cofibrant
homotopy group objects in $\C$, $\F$ a category of finite simplicial sets, and $\LL$ a
set of cones in $\C$. The associated \emph{dual enriched sketch}
\w{\TuA=\bTh\sp{(\uA,\K,\LL)}} is the smallest full sub-simplicial category of
$\C$ containing $\uA$ and closed under the operations \w{(-)\sp{K}} for \w{K\in\F}
and taking limits of the cones in $\LL$. In this setting:
\begin{enumerate}
\renewcommand{\labelenumi}{(\arabic{enumi})}
\item An \emph{$\uA$-dual presheaf} is a pointed simplicial functor
  \w[.]{\fX:\TuA\to\Sa}  The category of $\uA$-dual presheaves is denoted by
  \w[,]{\Sa\sp{\TuA}} and the value of $\fX$ at \w{\bA\in\TuA} will again be written
  \w[.]{\fX\lin{\bA}}
\item A \emph{\sAma} is a $\uA$-dual presheaf $\fX$ for which the natural maps
\begin{myeq}\label{eqloopprod}
\fX\lin{\bA\sp{K}}~\to~\fX\lin{\bA}\sp{K}~\hspace{9mm}\text{and}\hspace{9mm}
\fX\lin{\lim\sb{i\in I}\,\bA\sb{i}}~\to~\lim\sb{i\in I}\,\fX\lin{\bA\sb{i}}
\end{myeq}
\noindent are isomorphisms for all \w[,]{\bA, \bA\sb{i}\in\TuA} \w{K\in\F}, and diagrams
$I$ in $\LL$. The subcategory of \sAma[s] will be denoted by \w[.]{\sMuA}
\item A \emph{\wAma} is a $\uA$-dual presheaf $\fX$ which is weakly equivalent to
a \sAma[,] so in particular, the maps of \wref{eqloopprod} are weak equivalences
(see Remark \ref{rwsma} above). The subcategory of \wAma[s]  will be denoted
by \w[.]{\wMuA}
\end{enumerate}
\end{defn}

\begin{example}\label{egdmapalg}
The main example of an enriched dual sketch we consider here is the \emph{$\Omega$-spectrum}
case, where \w{\C=\Sa} and \w{\uA=\{\bA\sb{n}\}\sb{n=-\infty}\sp{\infty}}  are the spaces
of an $\Omega$-spectrum $A$ (in the sense of \cite{BFrieH}). The category
$\F$ then consists of the inclusions \w[,]{i\sb{0},i\sb{1}:\Delta[0]\hra\Delta[1]}
and the cone collection $\LL$ contains all products of cardinality $<\lambda$ for some fixed
limit cardinal $\lambda$ and the pullback squares
\mydiagram[\label{eqpathloop}]{
\ar @{} [drr] |(0.25){\framebox{\scriptsize{PB}}}
P\bA \ar@{^{(}->}[rr] \ar@{->>}[d]\sb{\simeq} &&
\bA\sp{\Delta[1]} \ar@{->>}[d]\sp{\simeq}\sb{\ev\sb{0}}&&
\ar @{} [drr] |(0.25){\framebox{\scriptsize{PB}}}
\Omega\bA \ar@{^{(}->}[rr]\sp{\iota\sb{\bA}} \ar@{->>}[d] &&
P\bA \ar@{->>}[d]\sb{\ev\sb{1}} \\
\ast \ar@{^{(}->}[rr] && \bA && \ast \ar@{^{(}->}[rr] && \bA
}
\noindent for any \w[.]{\bA\in\TuA}
Thus a \sAma $\fX$ will take the two pullback squares of \wref{eqpathloop}
to those of \wref[.]{eqcotolim}

More generally, one might take any set of $\Omega$-spectra \wh in particular,
the set of all $A$-module spectra of bounded cardinality, for a fixed ring spectrum $A$.
\end{example}

\begin{defn}\label{drma}
For any dual enriched sketch \w[,]{\TuA} the \emph{realizable} \sAma $\fX$
associated to \w{\bY\in\C}  has \w{\fX\lin{\bA}:=\map\sb{\C}(\bY,\bA)} for
each \w[.]{\bA\in\TsS} We will denote it by \w[.]{\fMA\bY}
When \w[,]{\bY\in\Obj\TuA} we again say that \w{\fMA\bY} is \emph{free}.
\end{defn}

The analogue of Lemma \ref{lfreema} also holds:

\begin{lemma}[cf.\ \protect{\cite[Lemma 1.12]{BSenM}}]\label{ldfreema}
If $\fY$ is an $\uA$-dual presheaf and \w{\fMA\bA} is a free \sAma (for \w[),]{\bA\in\TuA}
there is a natural isomorphism
$$
\Phi:\map\sb{\Sa\sp{\TuA}}(\fMA\bA,\fY)~\xra{\cong}~\fY\lin{\bA}~,
$$
\noindent with \w{\Phi(\fff)=\fff(\Id\sb{\bA})\in\fY\lin{\bA}\sb{0}}
for any \w[.]{\fff\in\Hom\sb{\Sa\sp{\TuA}}(\fMA\bA,\fY)}
\end{lemma}

\begin{defn}\label{ddpalg}
  As in \S \ref{dpalg}, given a dual enriched sketch \w[,]{\TuA} the corresponding
  ``algebraic'' sketch \w[,]{\ThuA:=\pi\sb{0}\TuA} whose models are now functors
  \w{\Lambda:\ThuA\to\Set} preserving all products among the cones listed in $\E$.
  These will be called \emph{\PAa[s]}, and their category will denoted by \w[.]{\PAAlg}
  Again, if $\fX$ is a (weak or strict) \ma[,] then \w{\pi\sb{0}\fX} is a \PAa[.]
  A \PAa is \emph{realizable} if it is isomorphic to \w{\pi\sb{0}\fMA\bY} for
  some \w[,]{\bY\in\C} and it is \emph{free} if it is of the form
  \w{\pi\sb{0}\fMA \bA} for \w[.]{\bA\in\Obj\TuA}
\end{defn}

\begin{example}\label{egasa}
  When \w{\A=\Fn} and \w[,]{\lambda=\aleph\sb{0}} \w{\TuA} is the simplicial category
  of finite type \ww{\Fp}-GEMs, and a \PAa is simply an unstable algebra over
  the mod $p$ Steenrod algebra (cf.\ \cite{SchwU}).
\end{example}

\begin{mysubsection}{Model categories of \ma[s]}
\label{smcma}
Like all categories of simplicial functors with small
indexing category, the (dual) presheaf categories \w{\Sa\sp{\TsS\op}} and
\w{\Sa\sp{\TuA}} have proper simplicial model category structures
(see \cite[13.1.14]{PHirM}), in which the fibrations and weak equivalences
are defined objectwise (see \cite[\S 1]{DKanSR}).
Thus a map \w{\fff:\fX\to\fY} of $\uS$-presheaves is a weak equivalence if for every
\w[,]{\bB\in\uS} \w{\fff\sb{\ast}:\fX\lin{\bB}\to\fY\lin{\bB}} is a weak equivalence
in $\C$ (as in \S \ref{dmapalg}).

By a suitable left Bousfield localization of \w{\Sa\sp{\TsS\op}}
and \w{\Sa\sp{\TuA}} we can obtain model categories for \wSma[s]
and \wAma[s] (i.e., model structures on the (dual) presheaf
category in which the latter are the fibrant objects). However,
since we cannot guarantee that these localized model structures
are right proper (cf.\ \cite[3.4.4]{PHirM}), they will not be used in this paper.
\end{mysubsection}

\begin{remark}\label{rfibrantma}
Note that since we assumed the objects of  \w{\TsS} are cofibrant, when $\bY$ is fibrant
the realizable $\uS$-presheaf \w{\fMS\bY} will be fibrant (that is, \w{\fMS\bY\lin{\bB}}
is a Kan complex for each \w[).]{\bB\in\TsS} Similarly, for $\A$-dual presheaves,
\w{\fMA\bY} is fibrant if $\bY$ is cofibrant in $\C$.
\end{remark}

\begin{mysubsection}{Model categories of simplicial \Pa[s]}
\label{smcspa}
Because both \PSa[s] (\S \ref{dpalg}) and \PAa[s] (\S
\ref{ddpalg}) are universal algebras in the sense of \cite[VI, \S
8]{MacLC} having an underlying graded group structure, there is a
model category structure on both the category \w{\PSAlg\sp{\Dop}}
of simplicial \PSa[s] and the category \w{\PAAlg\sp{\Dop}} of
simplicial \PAa[s]. In both cases  a map \w{f:\Ud\to\Vd} of
simplicial \Pa[s] is a weak equivalence (respectively,
fibration) if and only if the map
\w{f\sb{\ast}:\Ud\lin{\bB}\to\Wd\lin{\bB}} is a weak equivalence
(respectively, fibration) of simplicial groups for each
\w[.]{\bB\in\Obj\bTh} The cofibrant objects are retracts of free
simplicial objects.
\end{mysubsection}

\begin{mysubsection}{Truncating \ma[s]}
\label{struncma}
Fix \w[.]{n\geq 0} Given a $\uS$-presheaf \w[,]{\fX:\TsS\op\to\Sa} we may post-compose
$\fX$ with the $n$-th Postnikov section functor \w{\Po{n}:\Sa\to\Sn{n}} to obtain a
new $\uS$-presheaf \w[,]{\Po{n}\fX} which we now think of as a continuous functor on
\w{\Po{n}\TsS} \wh that is, the sketch enriched in \w{\Sn{n}} obtained from \w{\TsS}
by applying \w{\Po{n}} to each mapping space.

This is simplest to describe when $\fX$ is fibrant (cf.\ \S \ref{rfibrantma}),
since then we can use the \wwb{n+1}coskeleton functor \w{\csk{n+1}:\Sa\to\Sa}
(which strictly preserves products) as our model for \w[.]{\Po{n}} Note that the mapping
spaces of \w{\TsS} are always fibrant, since we assumed that all its objects are both
fibrant and cofibrant. In the general case, we must first apply a fibrant replacement
functor to $\fX$ in the model category \w{\Sa\sp{\TsS\op}} of \S \ref{smcma}.

The category of \emph{$n$-truncated $\uS$-presheaves} will be denoted by
\w[,]{\Sn{n}\sp{\TsS\op}\subset\Sa\sp{\TsS\op}} with the truncation functor
\w[.]{\gn{n}:\Sa\sp{\TsS\op}\to\Sn{n}\sp{\TsS\op}}

If $\fX$ is a (strict or weak) \Sma[,] this usually will not be true  of
\w[,]{\Po{n}\fX} since in general
\begin{myeq}\label{eqtruncma}
\Po{n}\map(\Sigma\bB,\bY)\simeq\Po{n}\Omega\map(\bB,\bY)\not\simeq
\Po{n-1}\Omega\map(\bB,\bY)\simeq\Omega\Po{n}\map(\bB,\bY).
\end{myeq}
\noindent Thus we must modify Definition \ref{dmapalg} as follows, assuming for
simplicity that the category $\F$ consists as above of the inclusions
\w[,]{i\sb{0},i\sb{1}:\Delta[0]\hra\Delta[1]} and the cocone collection $\E$ contains
all coproducts of cardinality $<\lambda$ for some fixed limit cardinal $\lambda$,
and the pushout squares \wref[:]{eqconesusp}

\begin{enumerate}
\renewcommand{\labelenumi}{(\arabic{enumi})}
\item An \emph{$n$-truncated \sSma} is an $n$-truncated $\uS$-presheaf $\fX$ for
which the natural maps of \wref{eqsuscoprod} are isomorphisms for
all \w[,]{\bB\in\TsS} \w[,]{K\in\F} and diagrams $I$ in $\E$, except for the right hand
square in \wref[,]{eqcotolim} where we have instead:
\begin{myeq}\label{eqtruncloop}
  \fX\lin{\Sigma\bB}~\to~\Po{n-1}\fX\lin{\Sigma\bB}~\xra{\cong}~\Po{n-1}\Omega\fX\lin{\bB}
  ~\xla{\cong}~\Omega\fX\lin{\bB}
\end{myeq}
\noindent where the first and last maps are the standard fibrations, the middle map
is the natural map of \wref[,]{eqsuscoprod} and \w{\Omega\fX\lin{\bB}} is an
\wwb{n-1}type by assumption, with the last map an isomorphism.

The full subcategory of $n$-truncated \sSma[s] will be denoted by \w[.]{\sMsSn{n}}
\item An \emph{$n$-truncated \wSma} is an $n$-truncated $\uS$-presheaf $\fX$
  weakly equivalent to an $n$-truncated \sSma[.] This implies that the maps of
  \wref[,]{eqsuscoprod} and the two right maps in \wref[,]{eqtruncloop} are weak
  equivalences (see Remark \ref{rwsma}).
The full subcategory of $n$-truncated \wSma[s] will be denoted by \w[.]{\wMsSn{n}}
\end{enumerate}
In particular, for any \w{\bY\in\C} we have the associated \emph{realizable}
$n$-truncated \sSma \w[,]{\Po{n}\fMS\bY} which is \emph{free} if
\w[,]{\bY\in\TsS} and the analogue of Lemma \ref{lfreema} still holds.
We define the $n$-truncated versions of $\uA$-dual presheaves and
(strict or weak) \Ama[s] dually.
\end{mysubsection}

%
%
\sect{Factoring functors through \ma[s]}
\label{cffma}

The first step in our program is to show that suitable homotopy functors \w{\bT:\C\to\D}
factor up to weak equivalence through an appropriate category of \ma[s:] in other words,
find an enriched sketch \w{\TsS}  and a functor \w[,]{\fT:\Sa\sp{\TsS\op}\to\D}
equipped with a natural weak equivalence \w[.]{\fT\circ\fMS\to \bT}
In fact, $\fT$ need not be defined on all of \w[;]{\Sa\sp{\TsS\op}}
it suffices if it is defined on the subcategory \w{\sMsS} of \sSma[s] where
\w{\fMS} takes values.

Dually, we could try to find a dual enriched sketch \w{\TuA}  and a functor
\w{\fT':\sMuA\to\D} with a natural weak equivalence \w[.]{\bT\to\fT'\circ\fMA}

\supsect{\protect{\ref{cffma}}.A}{Realizing mapping algebras}

The simplest way to define such a functor $\fT$ is in the case where every \sSma
$\fX$ is (functorially) realizable.  Essentially, the only case where this is known
to be true is when \w{\C=\Tz} and \w[.]{\uS=\{\bS{n}\}\sb{n=1}\sp{\infty}}
We briefly summarize the construction of \cite[\S 9]{BBlaC}  (based on that of
\cite[\S 2]{StoV}):

\begin{mysubsection}{The Stover construction}\label{sstovcons}
Recall that for a pointed Kan complex \w[,]{K\in\Sa}
the path space \w{PK} is given by \w[,]{(PK)\sb{n}:=\{x\in
K\sb{n+1}~:\ d\sb{1}\dotsc d\sb{n+1} x=\ast\}} with re-indexed
face and degeneracy maps, and the universal fibration \w{p:PK\to
K} is induced by \w{d\sb{0}} (cf.\ \cite{MoorS}). Thus when $K$ is
a simplicial group, the map on $0$-simplices
\w{p\sb{0}:(PK)\sb{0}\to K\sb{0}} suffices to compute
\w[.]{\pi\sb{0}K} We therefore choose the category \w{\G=\Grp\sp{\Dop}} of
simplicial groups as our model $\C$ for the homotopy theory of
pointed connected spaces, and set
\w{\uS:=\{\gS{n}\}\sb{n=1}\sp{\infty}} (where
\w[,]{\gS{n}:=F\bS{n-1}} as a free simplicial group, is a strict
cogroup object modelling the $n$-sphere in $\G$). For any limit
cardinal $\lambda$, the resulting enriched sketch
\w{\TsS=\TsS\sp{\lambda}} then has a strict mapping algebra
functor \w{\fMS:\G\to\sMsS} with each \w{\fMS\bY\lin{\bB}} a
simplicial group (though the structure maps are just maps of
pointed simplicial sets, in general).
\end{mysubsection}

\begin{defn}\label{dstover}
Let \w{\Gamma:=\bo\sp{\NN}} be the category consisting of a countable collection of
arrows, indexed by the objects of $\uS$, and \w{\Seta\sp{\Gamma}} the category
of $\Gamma$-indexed diagrams \w{\Phi:=(\phi\sb{n}:E\sb{n}\to F\sb{n})\sb{n\in\NN}}
in pointed sets, called \emph{arrow sets}.
We have a forgetful functor \w[,]{\rho:\sMsS\to\Seta\sp{\Gamma}} with
\w[.]{(\rho\fX)\sb{n}=(p\sb{0}:(P\fX\lin{\gS{n}}\sb{0}\to(\fX\lin{\gS{n}})\sb{0})}
In fact, \w{\rho\fX} is defined for any presheaf \w[,]{\fX:\TsS\op\to\Sa}  but we are only
interested in the composite \w[.]{\RS:=\rho\fMS:\G\to\Seta\sp{\Gamma}}
This has a left adjoint \w[,]{\LS:\Seta\sp{\Gamma}\to\G}
which assigns to an arrow set \w{\Phi=(\phi\sb{n}:E\sb{n}\to F\sb{n})\sb{n\in\NN}} the coproduct
\begin{myeq}\label{eqadjunction}
\LS \Phi~:=~\coprod\sb{n\in\NN}~\coprod\sb{f\in F\sb{n}}~~~\Ql{f}~,
\end{myeq}
\noindent where we define \w{\Ql{f}} for \w{f\in F\sb{n}} as follows:
\begin{enumerate}
\renewcommand{\labelenumi}{(\alph{enumi})~}
\item If \w[,]{\ast\neq f\in\Image\phi\sb{n}} then
\w{\Ql{f}} is defined by the pushout square
\mydiags[\label{eqsmallestpo}]{
  {\displaystyle\coprod\sb{e\in\phi\sb{n}\sp{-1}(f)}}~\gS{n}\sb{(e)} \ar[d]\sp(0.6){\coprod i\sb{(e)}}
  \ar[rr]\sp(0.6){\fold} && \gS{n}\sb{(f)} \ar[d]\\
  {\displaystyle\coprod\sb{e\in\phi\sb{n}\sp{-1}(f)}}~C\gS{n}\sb{(e)}\ar[rr] && \Ql{f}
}
\noindent in $\G$ (where \w{i:\gS{n}\to C\gS{n}} is the inclusion into the cone,
and $\fold$ is the fold map).
\item If \w[,]{f\not\in\Image \phi\sb{n}} we set \w[.]{\Ql{f}:=\gS{n}}
\item If \w[,]{f=\ast} we set
$$
  \Ql{f}~:=~\coprod\sb{\ast\neq e\in\phi\sb{n}\sp{-1}(\ast)}~\Sigma\gS{n}\sb{(e)}~.
$$
\end{enumerate}
Compare \cite[\S 2]{BSenM} and \cite[\S 2]{StoV}, where the comonad \w{\VS=\LS\RS:\G\to\G}
(or rather, its analogue for \w[)]{\Tz} was used to construct functorial resolutions
of pointed connected spaces by wedges of spheres.

Note that each \w[,]{\Ql{f}}  and thus \w[,]{\LS \Phi} is a strict cogroup object
in $\G$ (fibrant and cofibrant) of the homotopy type of a wedge of spheres.
If $\lambda$ is any limit cardinal, we define a \emph{$\lambda$-Stover space}
to be any pushout of the form \wref[,]{eqsmallestpo} with \w{\phi\sb{n}\sp{-1}(f)}
replaced by any set $T$ of cardinality $<\lambda$.  Let \w{\TsSt=\TsSt\sp{\lambda}}
denote  a skeleton of the sub-simplicial category of $\G$ whose objects are coproducts
of $\lambda$-Stover spaces over indexing sets of cardinality \www[.]{<\lambda}
This is an enriched sketch, with $\F$ as in \S \ref{egmapalg}, and $\E$ consisting
of the coproducts of cardinality $<\lambda$ in \w[,]{\TsSt} together with
the pushout squares of \wref{eqconesusp} and \wref[.]{eqsmallestpo} The category of the
corresponding \sMa[s,] called \emph{\sStma[s,]} will be denoted by \w[,]{\sMst} with
\w{\fMst:\G\to\sMst} the \emph{\sStma} functor.
\end{defn}

\begin{mysubsection}{The algebra structure}\label{salgst}
Since each sphere \w{\gS{n}\in\G} is in particular a Stover space, \w{\TsS=\TsS\sp{\lambda}}
is a full simplicial subcategory of \w[,]{\TsSt=\TsSt\sp{\lambda}} with \w{\iota:\TsS\hra\TsSt}
the inclusion, inducing the restriction \w[.]{\iota\sp{\ast}:\sMst\to\sMsS} Write
  \w{\wrh:\sMst\to\Seta\sp{\Gamma}} for the composite \w[.]{\rho\circ\iota\sp{\ast}}

  We claim that for every \sStma $\fX$, the arrow set \w{\rho\fX} has a natural
\ww{\TS}-algebra structure map \w{h:\TS\rho\fX\to\rho\fX} for the monad
\w{\TS=\RS\LS:\Seta\sp{\Gamma}\to\Seta\sp{\Gamma}} (see \cite[\S 4.1]{BorcH2}).
If we set \w{\K:=\LS\circ\wrh:\sMst\to\G} and \w[,]{\VS:=\fMst\circ\K} we may display the various
functors defined in the following commuting diagram:
\mydiagram[\label{eqfiunctors}]{
  &&&& \sMst \ar[dr]\sb{\iota\sp{\ast}} \ar@(ul,ur)^{\VS=\fMst\K}
  \ar@/^7.5pc/[dd]^{\wrh} \ar@/^5.5pc/[lllld]^{\K} \\
  \G \ar[rrrru]\sp{\fMst} \ar[rrrrr]\sp{\fMS} \ar@/^5.5pc/[rrrrd]^{\RS}
  &&&&& \sMsS \ar[dl]\sb{\rho}  \\
 &&&& \Seta\sp{\Gamma} \ar[llllu]\sb{\LS} \ar@(dr,dl)^{\TS=\RS\LS}
}
\noindent In this setting we have a stronger statement (cf.\ \cite[\ 9.19]{BBlaC}):
\end{mysubsection}

\begin{lemma}\label{lalgstr}
Every \sStma $\fX$ has a natural map \w{\xi\sb{\fX}:\VS\fX\to\fX} making the
following diagram
\mydiagram[\label{eqalgebra}]{
\VS\VS\fX \ar[rr]\sp(0.55){\xi\sb{\VS\fX}} \ar[d]\sb{\VS\xi\sb{\fX}} &&
\VS\fX \ar[d]\sp{\xi\sb{\fX}} \\
\VS\fX \ar[rr]\sb{\xi\sb{\fX}} && \fX
}
\noindent commute in \w[,]{\sMst} where \w{\xi\sb{\VS\fX}=\fMst\vare\sb{\K\fX}} for
\w{\vare:\K\fMst\to\Id} the counit of the comonad \w[.]{\LS\RS}
\end{lemma}

The structure map \w{h:\TS\wrh\fX\to\wrh\fX} is then given by \w[,]{\wrh(\xi\sb{\fX})} since
\w{\TS\circ\wrh=\wrh\circ\VS} (see \wref[).]{eqfiunctors}

\begin{proof}
Let \w{\D\sp{i}} denote either \w{\gS{n\sb{i}}} or \w{C\gS{n\sb{i}}} in $\G$.
\begin{enumerate}
\renewcommand{\labelenumi}{(\alph{enumi})~}
\item Recall that \w{\K\fX} is defined for any \sStma $\fX$ by the colimit
  \wref[,]{eqsmallestpo} which we may write as \w[,]{\colim\sb{i}\D\sp{i}\sb{f\sb{i}}}
  where \w[.]{f\sb{i}\in\fX\lin{\D\sp{i}}\sb{0}}
Since \w[,]{\K\fX\in\TsSt} the \sStma \w{\VS\fX} is free, so to define the algebra structure
map \w{\xi\sb{\fX}:\VS\fX\to\fX} we need only specify
\w[.]{\xi\sb{\fX}(\Id\sb{\K\fX})\in\fX\lin{\K\fX}\sb{0}} But $\fX$ takes the colimit of
\wref{eqsmallestpo} to a limit, so \w{\xi\sb{\fX}(\Id\sb{\K\fX})} is determined by the
elements \w[.]{f\sb{i}\in\fX\lin{\D\sp{i}}\sb{0}} We therefore write
\w[,]{\xi\sb{\fX}(\Id\sb{\K\fX})=\bigoplus\sb{i}\,f\sb{i}}
where $\oplus$ indicates that we are using the duality \wref{eqsuscoprod} between
the colimits and the limits.
\item Similarly, for any \w{Y\in\G} we have
  \w[.]{\K\fMst Y=\colim\sb{j}\D\sp{j}\sb{g\sb{j}}}
  The counit \w{\vare\sb{Y}:\K\fMst Y\to Y} is again determined by the
indexing maps as \w[,]{\vare\sb{Y}=\colim\sb{j} g\sb{j}}  with the induced map
\w{\fMst\vare\sb{Y}:\V\fMst Y\to\fMst Y} sending \w{\Id\sb{\K\fMst Y}} in
\w{\fMst\K\fMst Y\lin{\K\fMst Y}\sb{0}} to \w{[\colim\sb{j} g\sb{j}]} in
\w[.]{\fMst Y\lin{\K\fMst Y}\sb{0}}

Thus when \w[,]{\fX=\fMst Y} the map \w{\xi\sb{\fX}} sends \w{\Id\sb{\K\fMst Y}} to
\w{\vare\sb{Y}=\colim\sb{j}g\sb{j}} in
\w[.]{\fX\lin{\K\fMst Y}\sb{0}=\map(\K\fMst Y,\,Y)\sb{0}} This means that
\w[;]{\xi\sb{\fMst Y}=\fMst\vare\sb{Y}} in particular, the top horizontal map in
\wref{eqalgebra} is \w[.]{\xi\sb{\VS\fX}}
\item To evaluate the top right composite
  \w[,]{\varphi:=\xi\sb{\fX}\circ\fMst\vare\sb{\K\fX}:\VS\VS\fX\to\fX}
note that  \w{\VS\VS\fX} is free on \w[,]{\K\VS\fX} so we need only specify
\w{\varphi(\Id\sb{\K\VS\fX})} in \w[.]{\fX\lin{\K\VS\fX}} Since  \w{\xi\sb{\fX}}
is a map of \sStma[s,] it sends
\w{[\colim\sb{j} g\sb{j}]\in\VS\fX\lin{\K\VS\fX}\sb{0}} (for \w{Y:=\K\fX}
in (b) above) to
\begin{myeq}\label{eqcollim}
[\bot\sb{j} g\sb{j}]\sp{\ast}(\xi\sb{\fX}(\Id\sb{\K\fX})~=~
\top\sb{j} \,g\sb{j}\sp{\ast}(\bigoplus\sb{i}\,f\sb{i})\hsm
\text{in}\hsm \fX\lin{\K\VS\fX}\sb{0}~.
\end{myeq}
\item Since \w{\VS\VS\fX} is free, the map \w{\VS\xi\sb{\fX}~:~\VS\VS\fX\to\VS\fX}
  is determined by where it sends \w{\Id\sb{\K \VS\fX}} in
  \w[,]{\VS\fX\lin{\K \VS\fX}\sb{0}=\map(\K\VS\fX,\,\K\fX)\sb{0}} namely, to
  \w[.]{\K\xi\sb{\fX}:\K\VS\fX\to\K\fX} Since
\w{\K\VS\fX=\colim\sb{j}\D\sp{j}\sb{g\sb{j}}}
  where the colimit is over all maps \w[,]{g\sb{j}:\D\sp{j}\to\K\fX} we see
  from the description of \w{\xi\sb{\fX}} above (and the construction of $\K$)
  that \w{\K\xi\sb{\fX}} sends \w{\D\sp{j}\sb{g\sb{j}}} to the copy of
  \w{\D\sp{j}} in the colimit defining \w{\K\fX} indexed by
\begin{myeq}\label{eqeacheta}
\xi\sb{\fX}(g\sb{j})~=~\xi\sb{\fX}(g\sb{j}\sp{\ast}(\Id\sb{\K\fX}))~=~
g\sb{j}\sp{\ast}(\xi\sb{\fX}(\Id\sb{\K\fX}))~=~g\sb{j}\sp{\ast}(\bigoplus\sb{i}\,f\sb{i})~
\end{myeq}
\noindent in \w[,]{\fX\lin{\D\sp{j}}\sb{0}} where
\w{\xi\sb{\fX}(\Id\sb{\K\fX})=\bigoplus\sb{i}\,f\sb{i}}  by (a).

Thus the element \w{\xi\sb{\fX}(\VS\xi\sb{\fX}(\Id\sb{\K \VS\fX}))} in
\w{\fX\lin{\K\VS\fX}\sb{0}} is determined by the fact that $\fX$ takes the colimit
\w{\colim\sb{j}\D\sp{j}\sb{g\sb{j}}} defining \w{\K\VS\fX} to a limit,  namely:
\begin{myeq}\label{eqlimeta}
  \xi\sb{\fX}(\VS\xi\sb{\fX}(\Id\sb{\K \VS\fX}))~=~\xi\sb{\fX}(\K\xi\sb{\fX})~=~
  \xi\sb{\fX}(\bot\sb{j} g\sb{j})~=~\top\sb{j}\xi\sb{\fX}(g\sb{j})~=~
  \top\sb{j}g\sb{j}\sp{\ast}(\bigoplus\sb{i}\,f\sb{i})~.
\end{myeq}
\end{enumerate}
We see from \wref{eqcollim} and \wref{eqlimeta} that the two composites agree
on \w[,]{\Id\sb{\K\VS\fX}} so they are equal.
\end{proof}

\begin{mysubsection}{The resolution model category of simplicial presheaves}
\label{smcsps}
For any set \w{\uS\subset\C} as in \S \ref{dmapalg}, consider the category
\w{(\Sa\sp{\TsS\op})\sp{\Dop}=\Sa\sp{\TsS\op\times\Dop}} of
\emph{simplicial $\uS$-presheaves} \wh that is, simplicial objects in the category of
$\uS$-presheaves. As noted in \S \ref{smcma},
the $\uS$-presheaf category \w{\Sa\sp{\TsS\op}} has a proper simplicial model
category structure. Moreover, the objects of $\uS$ are homotopy cogroup objects in
$\C$, as are their colimits under $\E$ as in \S \ref{egmapalg}. Therefore, as in
\cite[\S 2]{JardB}, there is a resolution model category structure on
\w[,]{\Sa\sp{\TsS\op\times\Dop}} for which the projectives of
\w{\Sa\sp{\TsS\op}} are the free \sSma[s.] A map \w{\fff:\fVd\to\fWd} of simplicial
$\uS$-presheaves is a weak equivalence in this model category if and only if it is
an \emph{\ww{E\sp{2}}-equivalence} \wh that is, if for each \w{\bB\in\TsS}
and \w[,]{t,s,\geq 0} the map
\w{\fff\sb{\ast}:\pi\sb{t}\sp{h}\pi\sb{s}\sp{v}\fVd\lin{\bB}
  \to\pi\sb{t}\sp{h}\pi\sb{s}\sp{v}\fWd\lin{\bB}}
is an isomorphism (the terminology comes from the \Eut{2} of the homotopy
spectral sequence of a simplicial space \wh cf. \cite{DKStoE}).

Note that if a simplicial presheaf \w{\fWd} is cofibrant, each \w{\fW\sb{n}} is
weakly equivalent to a coproduct of free \sSma[s,] so in particular it is a \wSma[.]
Moreover, in order for \w{\fWd} to be a resolution of a \wSma $\fX$, in particular
\w{\pi\sb{0}\fWd} must be a resolution of \w{\pi\sb{0}\fX} in the model category
of simplicial \PSa[s]   (see \S \ref{smcspa}), so that the augmented simplicial group
\w{\pi\sb{0}\fWd\lin{\bB}\to\pi\sb{0}\fX\lin{\bB}} is weakly contractible for
any \w[.]{\bB\in\TsS}

We observe also that \w{\Sa\sp{\TsS\op\times\Dop}} has a Reedy model category structure,
with weak equivalences and fibrations defined at each simplicial space \w{\fVd\lin{\bB}}
for every \w{\bB\in\TsS} (see \cite[\S 15.3]{PHirM}).

Since \w{\Po{n}} is a nullification, \w{\Sn{n}} is still right proper (see
\cite[Theorem 9.9]{BousTH}, so we have an analogous resolution model category structure
on the category \w{\Sn{n}\sp{\TsS\op}} of $n$-truncated simplicial $\uS$-presheaves
(\S \ref{struncma}).
\end{mysubsection}

We deduce the following enhancement of \cite[Proposition 9.23]{BBlaC}:

\begin{thm}\label{trssma}
There is a \emph{realization functor} \w[,]{N:\sMst\to\G} equipped with natural
weak equivalences \w{\theta:N\circ\fMst\to\Id\sb{\G}} and
\w[.]{\zeta:\fMst\circ N\to\Id\sb{\sMst}}
\end{thm}

\begin{proof}
Given a \sStma \w[,]{\fX\in\sMst} iterating the comonad \w{\U:=\LS\RS:\G\to\G} on
\w{Y:=\K\fX=\LS\wrh\fX} yields an augmented simplicial space \w{\Zd\to Y} with
\w{\bZ\sb{n}:=\U\sp{n+1}Y} and \w{d\sb{i}:\bZ\sb{n}\to\bZ\sb{n-1}} given by as usual by
\w{\U\sp{i}\vare\sb{\U\sp{n-i}Y}} (cf.\ \cite[\S 8.6.4]{WeibHA}).

Since by \wref{eqfiunctors} \w{\U=\LS\RS=\K\fMst}  and
\w[,]{\VS=\fMst\K} we have a simplicial \sStma \w[,]{\fWd=\fMst\Zd} which augments
to $\fX$ via \w[,]{\xi\sb{\fX}:\fMst Y=\VS\fX\to\fX} by Lemma \ref{lalgstr}.
Applying $\K$ to \w{\fWd\to\fX} recovers \w[,]{\Zd\to Y} but now with an extra
degeneracy in each simplicial dimension coming from the unit \w{\eta:\Id\to\TS=\RS\LS}
of the corresponding monad, as well as an extra face map, obtained by iterating
$\U$ on \w[.]{\K\xi\sb{\fX}:\K\VS\fX=\bZ\sb{1}\to\K\fX=\bZ\sb{0}}
By commutativity of \wref[,]{eqalgebra} we see that \w{\Zd\to Y} is in fact
the d\'{e}calage of a simplicial space \w{\Xd} (see \cite{IlluC2}).
Moreover, applying \w{\fMst} to \w{\Xd} yields an augmented (free) simplicial
\sStma \w{\fMst\Xd\to\fX} which is a resolution of $\fX$ in the sense of \S \ref{smcsps}.

This shows that the Quillen-\-Bousfield-\-Fried\-lan\-der spectral sequence
for \w{\Xd} (see \cite{QuiS} and \cite[Theorem B.5]{BFrieH}) collapses, so that
\w{N\fX:=\|\Xd\|} realizes $\fX$ up to weak equivalence. Noting that \w{\Xd} is obtained
by applying $\K$ to \w[,]{\zeta\sb{0}:\fMst\Xd\to\fX} and that \w{\fMst\Xd}
is constructed by iterating $\VS$ on $\fX$ (together with \w[),]{\xi\sb{\fX}}
we have described a functorial procedure for realizing any \sStma $\fX$.
The natural weak equivalence $\zeta$ is induced by the augmentation
\w[,]{\zeta\sb{0}} while $\theta$ comes from the counit of the Stover comonad.
\end{proof}

\begin{cor}\label{crssma}
Any homotopy functor \w{\bT:\G\to\D} to a model category $\D$ induces a functor
  \w{\fT:=\bT\circ N:\sMst\to\D} equipped with a natural weak
equivalence \w[.]{\vartheta=\bT\theta:\fT\circ\fMst\to\bT}
\end{cor}

\supsect{\protect{\ref{cffma}}.B}{Realizing dual mapping algebras}

To dualize the results of \S \ref{cffma}.A, we want a setting where every
\sAma $\fX$ is functorially realizable. Again we have only one case where
this is known to be true, when \w{\C=\Sr} (or similar model categories
for pointed connected spaces) and $\uA$ consists of certain simplicial
$R$-modules for some commutative ring $R$.

\begin{defn}\label{detr}
In general, we must include in the corresponding enriched sketch \w{\TuA} all
$R$-module GEMs up to a certain cardinality. In particular, when \w{\C=\Sr} we
let \w[,]{\TuR=\TuR\sb{\lambda}:=\sMRla} be the full subsimplicial category of
$\C$ consisting of all simplicial $R$-modules of cardinality $<\lambda$,
for some limit cardinal $\lambda$ (determined as in \cite[\S 3.B]{BSenM}).
The corresponding dual \ma[s] will be
called \emph{\sRma[s]} (or \emph{\Rma[s,]} for short), and the category of such
will be denoted by \w[,]{\sMuR} with \w{\fMR:\C\op\to\sMuR} the \emph{realizable} \Rma[s.]
\end{defn}

\begin{mysubsection}{The dual Stover construction}\label{sdstovcons}
As in \S \ref{sstovcons}, we have a forgetful functor
\w[,]{\rho:\sMuR\to(\Seta\sp{\Gamma})\op} with
\w[.]{(\rho\fX)\sb{n}=(p\sb{0}:(P\fX\lin{\KR{n}}\sb{0}\to(\fX\lin{\KR{n}})\sb{0})}
The composite \w{\LR:=\rho\fMR:\C\to(\Seta\sp{\Gamma})\op}  has a right adjoint
\w[,]{\RR:(\Seta\sp{\Gamma})\op\to\C} with
\w{\RR\Phi:=\prod\sb{n\in\NN}~\prod\sb{f\in F\sb{n}}~\Qu{f}} for
any arrow set \w{\Phi=(\phi\sb{n}:E\sb{n}\to F\sb{n})\sb{n\in\NN}}

When $R$ is a field, we define \w{\Qu{f}} for \w{f\in F\sb{n}} by the pullback square
\mydiags[\label{eqsmallestpb}]{
\Qu{f} \ar[rr] \ar[d] &&
{\displaystyle\prod\sb{\phi\sb{n}\sp{-1}(f)}}~P\KR{n} \ar[d]\sp(0.6){\prod p\sb{\KR{n}}}\\
\KR{n}\ar[rr]^(0.4){\diag} &&
{\displaystyle\prod\sb{\phi\sb{n}\sp{-1}(f)}}~\KR{n}
}
\noindent if \w[,]{\ast\neq f\in\Image\phi\sb{n}}
while \w{\Qu{f}:=\KR{n}} if \w[.]{f\not\in\Image \phi\sb{n}} If \w[,]{\phi=\ast}  we set
\w{\Qu{f}~:=~\prod\sb{\phi\sb{n}\sp{-1}(\ast)\setminus\{\ast\}}~\Omega \KR{n}}
(compare \wref[).]{eqsmallestpo}

Again, for any limit cardinal $\Lambda$ we define a \emph{$\lambda$-$R$-Stover space}
to be any pullback of the form \wref[,]{eqsmallestpb} with \w{\phi\sb{n}\sp{-1}(f)}
replaced by any set $T$ of cardinality $<\lambda$.
When $R$ is not a field, we need to use the more complicated \emph{modified Stover
  construction} of \cite[\S 3.A]{BSenM} instead of the above.

We denote by \w{\TuSRl} the corresponding dual enriched sketch,
with $\F$ as in \S \ref{egdmapalg}, and $\LL$  consisting of products of
cardinality $<\lambda$ in \w[,]{\TuSRl} together with the pullback squares
of \wref{eqpathloop} and \wref[.]{eqsmallestpb} The category of the corresponding
dual \sMa[s,] called \emph{\dsStma[s,]} will be denoted by
\w[,]{\dsMst} with \w{\fMRst:\C\to\dsMst} the \emph{\dsStma} functor.

Since each \w{\KR{n}} is in particular an $R$-Stover space,  \w{\TuR\sb{\lambda}}
is a full simplicial subcategory of \w[,]{\TuSRl} with
\w{\iota:\TuR\sb{\lambda}\hra\TuSRl} the inclusion, inducing the
restriction \w{\iota\sp{\ast}:\dsMst\to\sMuR} as in \S \ref{salgst}.
Writing \w[,]{\VR:=\fMRst\circ\RR\circ\rho\circ\iota\sp{\ast}:\dsMst\to\dsMst}
we obtain the following categorical dual of Lemma \ref{lalgstr}
(compare \cite[Proposition 2.19]{BSenM}):
\end{mysubsection}

\begin{lemma}\label{ldalgstr}
Every \dsStma $\fX$ has a natural map \w{\zeta\sb{\fX}:\VR\fX\to\fX} making the
following diagram commute in \w[:]{\dsMst}
\mydiagr[\label{eqdalgebra}]{
\VR\VR\fX \ar[rr]\sp(0.55){\zeta\sb{\VR\fX}} \ar[d]\sb{\VR\zeta\sb{\fX}} &&
\VR\fX \ar[d]\sp{\zeta\sb{\fX}} \\
\VR\fX \ar[rr]\sb{\zeta\sb{\fX}} && \fX
}
\end{lemma}

\begin{defn}\label{dewr}
For any commutative ring $R$, we denote by \w{\SR} the full subcategory of
$R$-good spaces in \w{\Sa} (cf.\ \cite[I, \S 5.1]{BKanS}), and by \w{\dsMsr} the full
subcategory of \w{\dsMst} consisting of those  \dsStma[s] which
are weakly equivalent to \w{\fMRst\bY} for some \w[.]{\bY\in\SR}
These will be called \emph{weakly $R$-good \dsStma[s.]}
\end{defn}

\begin{remark}\label{rmcspa}
By \cite[\S 15.3]{PHirM}), \w{\Sa\sp{\TuA\times\Dop}} and
\w{\Sn{n}\sp{\TuA\times\Dop}} have Reedy model category structures, with
weak equivalences and cofibrations defined at each simplicial space \w{\fWd\lin{\bB}}
for each \w[.]{\bB\in\TuA}

As in \S \ref{smcsps}, there is also a resolution model category structure on the category
\w{(\Sa\sp{\TuA})\sp{\Dop}=\Sa\sp{\TuA\times\Dop}} of simplicial
dual $\uA$-presheaves. Again, if a simplicial presheaf \w{\fWd} is cofibrant,
each \w{\fW\sb{n}} is weakly equivalent to a coproduct of free \sAma[s,] so it is
a \wAma[,] and \w{\fWd\to\fX} is a resolution of \wAma[s] only if
\w{\pi\sb{0}\fWd\to\pi\sb{0}\fX} is a resolution of  \PAa[s.]

Since \w{\Sn{n}} is still proper, we also have a resolution model category
structure on the category \w{\Sn{n}\sp{\TuA\times\Dop}}
of $n$-truncated simplicial dual $\uA$-presheaves (\S \ref{struncma}).
\end{remark}

The Eckmann-Hilton dual of Theorem \ref{trssma} has the following more involved form:

\begin{thm}\label{tdrssma}
Let $R$ be any commutative ring, \w[,]{\C=\Sa}  and $\fX$ a \sRma (for \w{\TuR=\sMRla}
as in \S \ref{detr}), which we assume to be a \dsStma[.]
\begin{enumerate}
\renewcommand{\labelenumi}{(\alph{enumi})}
\item There is a functor associating to $\fX$ a cosimplicial object
\w{\Wu\in\Sa\sp{\Delta}} with each \w{W\sp{n}} in \w[,]{\sMRla} equipped with
a natural augmentation of \Rma[s] \w[,]{\vare:\fMR\Wu\to\fX} such that
\w{\pi\sb{0}\fMR\Wu\to\pi\sb{0}\fX\lin{M}} is a simplicial resolution of \PAa[s.]
\item If \w{\fX\in\dsMsr} is weakly equivalent to \w{\fMRst\bY}
  (for some $R$-good space $\bY$), then \w{\Tot\Wu} is homotopy equivalent to
  the $R$-completion of $\bY$ (so in particular \w{\Tot\Wu} realizes $\fX$
  up to weak equivalence).
\item When $R$ is a field, we can start with any \sAma $\hfX$ (for
  \w{\uA=\Rn} in \S \ref{egdmapalg}).  If it extends to a \dsStma $\fX$ as defined
  in \S \ref{sdstovcons}, and then (a) and (b) hold.
\item When \w{R=\Fp} or $\QQ$, and $\fX$ is simply connected
  (that is, letting \w{\uA=\{\KR{n}\}\sb{n=2}\sp{\infty}} in \S \ref{egdmapalg}),
  \emph{any} \Rma (for a suitable limit cardinal $\lambda$) is weakly equivalent
  to \w{\fMRst\bY} for some simply connected $\bY$, unique up to $R$-equivalence.
\end{enumerate}
\end{thm}

\begin{proof}
This follows from various results in \cite{BSenM}:
\begin{enumerate}
\renewcommand{\labelenumi}{(\alph{enumi})}
\item This is \cite[Proposition 3.9]{BSenM}.
\item This is \cite[Theorem 3.26]{BSenM}.
\item This combines \cite[Proposition 2.23]{BSenM} and \cite[Theorem 2.30]{BSenM},
  using the fact that a weak equivalence of \dsStma[s] \w{\fff:\fX\to\fY}
  induces weak equivalence (in the model category of \cite[\S 3]{BousCR})
  between the corresponding cosimplicial spaces (see \cite[\S 7.7]{BousCR}).
\item This is \cite[Theorem 4.23]{BSenM} (when \w[)]{\lambda=\aleph\sb{0}} or
  \cite[Theorem 4.28]{BSenM} (otherwise).
\end{enumerate}
\end{proof}

\begin{cor}\label{cdrssma}
If $R$ is any commutative ring, there is a \emph{realization} functor
\w{N:(\dsMsr)\op\to\Sa} with a natural weak equivalence \w[.]{\vare:\Id\to N\circ\fMRst}
Thus any functor \w{\bT:\SR\to\D} (see \S \ref{dewr}) to a model category
$\D$ which preserves $R$-equivalences induces a functor
\w{\fT:=\bT\circ N:(\dsMsr)\op\to\D} equipped with a natural weak equivalence
\w[.]{\vartheta=\bT\vare:\fT\to\fT\circ\fMRst}
\end{cor}

\begin{proof}
We set \w[,]{N:=\Tot\Wu}  where \w{\fX\mapsto\Wu} is the functor
of Theorem \ref{tdrssma}. Once we know that $\fX$ is weakly
$R$-good (see \S \ref{dewr}), the natural augmentation
\w{\vare:\fMRst N\to\fX} is a weak equivalence by Theorem
\ref{tdrssma}(b) or (c).
\end{proof}

\begin{example}\label{egdrssma}
For any \w{\bZ\in\Sa} with \w{\bT:\SR\to\Sa} the functor
\w[,]{\mapa(\bZ,-)} the induced functor \w{\fT:=\bT\circ N:\dsMsr\to\Sa}
has the property that if \w{\fZ:=\fMRst\bZ} and $\fX$ is the realizable \dsStma
\w{\fMRst\bY} for some $R$-good space $\bY$, then \w{\fT(\fX)} is weakly equivalent to
\w[.]{\fZ\lin{\bY}}

Thus the $n$-truncation \w{\Po{n}\fT} (cf.\ \S \ref{struncma}) when
evaluated at \w[,]{\fX=\fMRst\bY} is determined by  \w[.]{\Po{n}\fZ} Moreover,
from the alternative description in \S \ref{raltma} we see that if
\w[,]{\bY\in\TuSRl} then \w{\fT(\fX)} corresponds to the $n$-truncated simplicial
category \w[,]{\Po{n}\cX} so that in fact \w[,]{\Po{n}\fT} when evaluated at
free \dsStma[s,] factors through the $n$-truncation.
\end{example}

%
%
\sect{Relative derived functors}
\label{chotdf}

Let \w{\bT:\C\to\D} be a homotopy functor between model categories of spaces.
We want to study the homotopy spectral sequence for the (co)simplicial object
obtained by applying $\bT$ to a (co)simplicial resolution of a space
\w[,]{\bY\in\C} using a relative version of the total
derived functor of the associated functor of \ma[s] $\fT$.

\begin{mysubsection}{Relative left and right derived functors}
\label{srdf}
If \w{T:\D\to\E} is a functor between model categories which
preserves weak equivalences of cofibrant objects, recall that
Quillen constructs the \emph{total left derived functor} \w{\bLL
T:\ho\D\to\ho\E} on an object \w{x\in\D} by applying $T$ to any
cofibrant replacement of $x$ (see \cite[I, \S 4]{QuiH}). In order
for this to work, $T$ need only be defined on the full subcategory
\w{\D\sb{\cf}} of all cofibrant objects in $\D$. In the spirit of
the Eilenberg-Moore ``relative homological algebra'' (see
\cite{EMoorF}), one could require only that $T$ be defined on some
full subcategory  $\PP$ of \emph{special} cofibrant objects in
\w{\D\sb{\cf}} (e.g., free, rather than projective, resolutions)
\wh as long as every object of $\D$ is weakly equivalent to an
object of $\PP$ (and $T$ still takes weakly equivalent objects of
$\PP$ to weakly equivalent objects in $\D$). Moreover, if are only
given a full subcategory \w{\D\sb{\PP}} of $\D$, closed under weak
equivalences, and every object of \w{\D\sb{\PP}} is weakly
equivalent to one in $\PP$, we still have \w[.]{\bLL
T:\ho\D\sb{\PP}\to\ho\E} Finally, $\E$ need not be a model
category \wh all we need is the localization
\w[,]{\gamma:\E\to\ho\E} with \w{\gamma\circ T} taking weak
equivalences to isomorphisms.

However, we shall be interested in a situation where we have two model category
structures on $\D$ \wh or perhaps only a subcategory \w{\Wc} of the given weak equivalences
$\cW$. This commonly occurs when our model category \w{(\D,\cW,\D\sb{\cf},\D\sb{\fib})}
is obtained by localizing another.

In this case, we shall assume that $\PP$ and \w{\D\sb{\PP}} satisfy the stronger
requirement that for each \w{x\in\D\sb{\cf}\cap\D\sb{\PP}} there is a map
\w{f:y\to x} in \w{\Wc} with \w[.]{y\in\PP}
If \w{T:\PP\to\E} is then a functor which preserves $\cW$-weak equivalences,
the \emph{relative left derived functor} of $T$ (with
respect to $\PP$ and \w[)]{\Wc} is the functor \w{\bLr T:\ho\D\to\ho\E}
defined on \w{z\in\D\sb{\PP}} by applying $T$ to $y$, where \w{g:x\to z}
is a cofibrant replacement (with respect to $\cW$) and \w{f:y\to x} in \w{\Wc}
is as above.

Dually, if we have full subcategories $\F$ of \w{\D\sb{\fib}} (the fibrant objects)
and \w{\D\sb{\F}} of $\D$, both closed under weak equivalences, and
\w[,]{\Wc\subseteq\cW} with the corresponding dual properties with respect
to a homotopy functor \w[,]{T:\F\to\E} the \emph{relative right derived functor}
\w{\bRr T:\ho\D\sb{\F}\to\ho\E} is defined analogously.
\end{mysubsection}

\begin{remark}\label{rappldf}
In the applications we have in mind, $\D$ will be a resolution model
category of simplicial \ma[s,] so the weak equivalences
$\cW$ in $\D$ are \ww{E\sp{2}}-equivalences. However, we also have a Reedy model
structure on $\D$, and the special weak equivalences \w{\cW'} will be
the \ww{E\sp{1}}-equivalences. The ability to apply the functor $T$ to a resolution
which is \ww{\cW'}-equivalent to \emph{any} cofibrant replacement (that is,
simplicial resolution) $y$ of an object \w{z\in\D\sb{\F}} provides the flexibility
we want in using particular resolutions \wh e.g., minimal \wh to calculate
\w[,]{(\bLr T)z} and eventually, the appropriate terms of the
spectral sequence.
\end{remark}

\supsect{\protect{\ref{chotdf}}.A}{Relative left derived functors
of \ma[s]}

For \w{\C=\G} and \w{\TsS} as in \S \ref{egmapalg}, let \w{\Wd} be a resolution
of \w{\bY\in\G} in the resolution model category structure on \w[.]{\G\sp{\Dop}}
Given a homotopy functor \w{\bT:\G\to\D} for $\D$ a ``category of spaces'' such
as \w[,]{\Tz} \w[,]{\Sa} or $\G$, we wish to study the homotopy spectral sequence
for the simplicial space \w{\bT\Wd\in\D\sp{\Dop}}

By applying the functor \w{\fMst:\C\to\sMst} of \S \ref{sdstovcons} to \w[,]{\Wd}
we obtain a simplicial \sStma \w{\fWd:=\fMst\Wd} which is a cofibrant replacement for
\w{\fX:=\fMst\bY} in the resolution model category structure on
\w{\Sa\sp{\TsS\op\times\Dop}} associated to the free \dsStma[s]
\w[.]{\{\fMst\bS{i}\}\sb{i=1}\sp{\infty}} By Corollary \ref{crssma},
there is functor \w[,]{\fT=\bT N:\sMst\to\D}  with a natural Reedy (that is,
levelwise) weak equivalence of simplicial spaces \w[.]{\vartheta:\fT\fWd\to\bT\Wd}

We want to calculate the total left derived functor of $\fT$ evaluated at $\fX$ by
applying $\fT$ to any resolution \w[.]{\fVd\to\fX} However, such an \w{\fVd}
is just a simplicial $\uS$-presheaf, and the functor $\fT$ is
only defined for \emph{strict} \Stma[s.] As explained in \S \ref{srdf}, our
solution to this difficulty is to show that any such \w{\fVd} is in fact
\ww{E\sp{1}}-equivalent to a simplicial \sSma \w[.]{\fWd}
For this purpose we require some additional notions from \cite[\S 1]{BJTurHI}:

\begin{mysubsection}{CW resolutions}
\label{scompres}
If $\E$ is any pointed complete category, the $n$-th \emph{Moore chains} object
of \w{\Gd\in\E\sp{\Dop}} is
\w[.]{C\sb{n}\Gd~:=~\cap\sb{i=1}\sp{n}\Ker\{d\sb{i}:G\sb{n}\to G\sb{n-1}\}}
The differential is
\w{\partial\sb{n}:=d\sb{0}\rest{\C\sb{n}\Gd}:C\sb{n}\Gd\to C\sb{n-1}\Gd}
and the \emph{cycles} objects is \w[,]{Z\sb{n}\Gd:=\Ker(\partial\sb{n})} with
\w{v\sb{n}:Z\sb{n}\Gd\to C\sb{n}\Gd} the inclusion.
These are defined for any restricted simplicial object
\w{\Gd\in\E\sp{\Dop\sb{+}}} (see \S \ref{snac}).

The $n$-th \emph{latching object} for \w{\Gd} is the colimit
\begin{myeq}\label{eqlatch}
L\sb{n}\Gd~:=~\colimit{\theta\op:\bbk\to\bbn}\,G\sb{k}~,
\end{myeq}
\noindent where $\theta$ ranges over the surjective maps \w{\bbn\to\bbk} in
$\Delta$ for \w[.]{k < n}

A simplicial object \w{\Gd\in\E\sp{\Dop}} is called a
\emph{CW object} if it is equipped with a \emph{CW basis}
\w{(\oG{n})\sb{n=0}\sp{\infty}} in $\E$ such that
\w[,]{G\sb{n}=\oG{n}\amalg L\sb{n}\Gd} and \w{d\sb{i}\rest{\oG{n}}=0}
for \w[.]{1\leq i\leq n} The $n$-th \emph{attaching map} for \w{\Gd} is defined
to be \w{\odz{G}{n}:=d\sb{0}\rest{\oG{n}}:\oG{n}\to C\sb{n-1}\Gd}
(which actually lands in \w[).]{Z\sb{n-1}\Gd}

When $\E$ is a suitable category of universal algebras, such as \w{\PSAlg}
(cf.\ \S \ref{dpalg}), a simplicial object
\w{\Vd\in\E\sp{\Dop}} with an augmentation to \w{\Lambda\in\C} is called a
\emph{CW resolution} if \w{\Vd\to\Lambda} is acyclic, with a CW basis
\w{(\oV{n})\sb{n=0}\sp{\infty}} having each \w{\oV{n}} free. Moreover, in this case
\w{\odz{V}{n}} surjects onto \w{Z\sb{n-1}\Vd} (where \w[).]{Z\sb{-1}\Vd:=\Lambda}

For \w[,]{\uS=\{\bS{i}\}\sb{i=1}\sp{\infty}} by \cite[Lemma 1.38]{BJTurHI}
every free simplicial \PSa (\S \ref{dpalg}) has a free CW basis. Moreover,
by \cite[Theorem 2.29]{BJTurHI}, every CW resolution \w{\Vd} of a realizable
\PSa \w{\Lambda=\pis\bY=\pi\sb{0}\fMS\bY} can be realized by an augmented
simplicial space \w[.]{\Wd\to\bY} Therefore, every free simplicial \PSa
resolution \w{\Vd\to\pis\bY} can be realized (non-canonically) by a
simplicial resolution of \sSma[s] \w[,]{\fWd\to\fMS\bY}
with \w[.]{\pi\sb{0}\fWd\cong\Vd} In order to apply the ideas of \S \ref{srdf},
we must show that any simplicial $\uS$-presheaf resolution \w{\fVd} of
\w{\fX=\fMS\bY} is Reedy weakly equivalent to a \sStma resolution \w[.]{\fWd}
To do so, we recall the following constructions from \cite{BJTurHI}:
\end{mysubsection}

\begin{mysubsection}{Sequential realizations}\label{sseqreal}
Assume given an enriched sketch \w{\TsS=\bTh\sb{(\uS,\F,\E)}}
in a pointed simplicial model category $\C$, as in \S \ref{dmapalg}, and
a CW-resolution \w{\Vd} of a realizable \PSa \w[,]{\Lambda=\piS\bY} with
CW basis \w[.]{\{\oV{n}\}\sb{n=0}\sp{\infty}} We define a
\emph{sequential realization of \w{\Vd}} (for $\bY$) to be a sequence $\cW$ of maps
\begin{myeq}\label{eqtower}
\WW{0}~\xra{\irn{0}}~\WW{1}~\xra{\irn{1}}~\WW{2}~\to~\dotsc~
\WW{n}~\xra{\irn{n}}~\WW{n+1}~\to~\dotsc~
\end{myeq}
\noindent between Reedy fibrant and cofibrant objects in
\w[,]{C\sp{\Dop}} such that for each \w[:]{n\geq 0}

\begin{enumerate}
\renewcommand{\labelenumi}{(\roman{enumi})~}
\item \w{\oW{n}\in\TsS} realizes the given CW basis \PAa \w[.]{\oV{n}}
\item There is an $n$-skeletal restricted simplicial object \w{\tWd{n}} with
\begin{myeq}\label{eqrsocp}
\tWn{k}{n}~=~\Wn{k}{n-1}\amalg\CsW{n-k-1}{n}\hsp \text{for}\hsm 0\leq k\leq n,
\end{myeq}
\noindent where by convention \w[,]{\CsW{0}{n}:=C\oW{n}}   \w[.]{\CsW{-1}{n}:=\oW{n}}
and \w[.]{\WW{-1}=\ast}
\item The face map \w{d\sb{0}\rest{\CsW{n-k-1}{n}}} is the map \w{F\sb{k}}
  in the commuting diagram
\mydiagram[\label{eqakfk}]{
  \Sigma\sp{n-k-1}\oW{n} \ar@{^{(}->}[rr]^{i\sp{k}} \ar[d]^{a\sb{k}} &&
  \CsW{n-k-1}{n} \ar@{->>}[rr]^{q\sp{k}} \ar[d]^{F\sb{k}} &&
  \Sigma\sp{n-k}\oW{n} \ar[d]^{a\sb{k-1}} \\
  Z\sb{k-1}\WW{n-1} \ar@{^{(}->}[rr]^{v\sb{k-1}} &&
  C\sb{k-1}\WW{n-1} \ar@{->>}[rr]^{\odz{}{k-1}} && Z\sb{k-2}\WW{n-1}
}
\noindent in which the top row is a strict cofibration sequence and the bottom
row a strict fibration sequence in $\C$. Thus \w{F\sb{k}} is a nullhomotopy
for \w[,]{v\sb{k-1}\circ a\sb{k}} which in turn defines \w[,]{a\sb{k-1}} using
\wref[.]{eqakfk} The first face map \w{d\sb{1}\rest{\CsW{n-k-1}{n}}} is the composite
\w[,]{\CsW{n-k-1}{n}\xra{q\sp{k}}\Sigma\sp{n-k}\oW{n}\xra{i\sp{k-1}}\CsW{n-k}{n}}
and \w{d\sb{i}\rest{\CsW{n-k-1}{n}}=0} for \w[.]{i>1}

We start with \w{F\sb{n}:\oW{n}\to C\sb{n-1}\WW{n-1}} a realization of the $n$-th
attaching map \w{\odz{V}{n}:\oV{n}\to C\sb{n-1}\Vd} for the given CW resolution, and
\w{a\sb{n-1}:=\odz{}{n-1}\circ F\sb{n}:\oW{n}\to Z\sb{n-2}\WW{n-1}} (with
\w{v\sb{n-2}\circ a\sb{n-1}} indeed nullhomotopic).
\item Let \w{\vWu{n}} be the pushout of the obvious maps
\begin{myeq}\label{eqrsopo}
\WW{n-1}~\leftarrow~\LL i\sp{\ast}\WW{n-1}~\to~\LL\tWd{n}~,
\end{myeq}
\noindent where \w{\LL:\C\sp{\Dzop}\to \C\sp{\Dop}} is the left adjoint of
\w[,]{i\sp{\ast}:\C\sp{\Dop}\to \C\sp{\Dzop}} as in \S \ref{snac}.
We then let \w{\WW{n}} be a Reedy fibrant and cofibrant replacement for \w[.]{\vWu{n}}
\item There is an augmentation \w{\bve{n}:\WW{n}\to\bY} realizing
\w{\Vd\to\Lambda} through simplicial dimension $n$ \wh that is, the $n$-truncation
of the augmented simplicial \PAa \w{\piA\WW{n}\to\piA\bY} is isomorphic to
the $n$-truncation of \w[.]{\Vd\to\Lambda}
\item The maps \w{\irn{n}} restrict to a trivial cofibration
\w{\irnk{n}{k}:\Wn{k}{n-1}\xra{\simeq}\Wn{k}{n}} for each \w[.]{0\leq k<n}
\end{enumerate}

It follows that \w{\Wd:=\colim\sb{n}\WW{n}\xra{\bv}\bY} is a simplicial resolution
in the resolution model category \w[.]{\C\sp{\Dop}}
See \cite[\S 2]{BJTurHI} for further details.
\end{mysubsection}

\begin{thm}\label{ttlwe}
For an enriched sketch \w{\TsS} as in \S \ref{dmapalg}, \w{\bY\in\C} fibrant, and
\w[,]{\fX:=\fMS\bY} let \w{\bet:\fVd\to\co{\fX}} be a trivial fibration with
\w{\fVd} cofibrant in \w[.]{\Sa\sp{\TsS\op\times\Dop}}
Then for any sequential realization $\cW$ of the \PSa resolution
\w{\pi\sb{0}\fVd\to\piS\bY} as in \S \ref{sseqreal}, there is
a Reedy weak equivalence of simplicial \wSma[s] \w[.]{\fff:\fMS\Wd\to\fVd}
\end{thm}

\begin{proof}
By \S \ref{smcsps},  the simplicial \PSa \w{\Vd:=\pi\sb{o}\fVd} is a free resolution
of the \PSa \w[,]{\Lambda:=\pi\sb{0}\fX} so it has a CW basis
\w{\{\oV{n}\}\sb{n=0}\sp{\infty}} by \cite[Lemma 1.38]{BJTurHI}, with
\w{\oV{n}=\pi\sb{0}\fMS\oW{n}} for some \w[.]{\oW{n}\in\Obj\TsS}
We may assume \w{\TsS} contains all simplicial groups
of the homotopy type of a (possibly trivial) wedge of objects of $\uS$ of
cardinality \www[.]{<\lambda}  This will ensure that all objects \w[,]{\Wn{k}{n}}
\w[,]{\tWn{k}{n}} \w[,]{\vWn{k}{n}} and so on, in \S \ref{sseqreal}
are in \w[.]{\TsS}

We construct $\fff$ by a double induction: in the outer induction, we construct
maps of simplicial \wSma[s] \w[.]{\fff\bp{n}:\fMS\WW{n}\to\fVd}
Assuming we have defined \w[,]{\fff\bp{n-1}} we need to extend it to a map  of
$n$-truncated restricted simplicial objects \w[,]{\tff{n}:\fMS\tWd{n}\to\fVd}
which we do by an inner downward induction on \w[.]{0\leq k\leq n}
Using Lemma \ref{lfreema}, we see from \wref{eqrsocp} that \w{\tff{n}\sb{k}} is
determined by an element \w{\of{k}{n}\in\fV\sb{k}\lin{\CsW{n-k-1}{n}}\sb{0}}
with \w{d\sb{i}\of{k}{n}=0} for \w[\vsm.]{i\geq 2}

\noindent\textbf{Step A.}\
To start the outer induction, note that since \w[,]{\Wn{0}{0}=\oW{0}} by
Lemma \ref{lfreema} the augmentation \w{\bve{0}:\fMS\WW{0}\to\fX} is determined by
an element \w[.]{e\in\fX\lin{\oW{0}}\sb{0}=\Hom(\oW{0},\bY)}
 Since \w{\bet:\fVd\to\co\fX} is a Reedy fibration (see \cite[\S 2]{JardB}),
 \w{(\bet\sb{0})\sb{\ast}:\fV\sb{0}\lin{\oW{0}}\to\fX\lin{\oW{0}}} is a fibration and
 in particular a surjection in \w[.]{\Sa} Moreover,
\w{\pi\sb{0}\fV\sb{0}\cong\piS\oW{0}}
 is a free \PSa[,] by our assumption on \w[,]{\fVd} so we have an element
 \w{\of{0}{0}\in\fV\sb{0}\lin{\oW{0}}\sb{0}} representing
 \w{\Id\in\pi\sb{0}\fV\sb{0}\lin{\oW{0}}} with \w{(\bet\sb{0})\sb{\ast}\of{0}{0}=e}
 by \cite[Lemma 15.9]{BJTurR}, as required\vsm.

\noindent\textbf{Step B.}\ Given \w[,]{\fff\bp{n-1}:\fMS\WW{n-1}\to\fVd} consider
the augmented simplicial space \w[:]{\Xd:=\fVd\lin{\oW{n}}\to\fX\lin{\oW{n}}} we
think of this as a bisimplicial set with vertical direction internal to each
\w[,]{\fV\sb{k}\lin{\oW{n}}\in\Sa} and horizontal direction corresponding to the
original simplicial direction of \w[.]{\fVd} The (split) inclusion
\w{j\sb{n}:\oV{n}\hra V\sb{n}} for the CW basis \PSa \w{\oV{n}=\pi\sb{0}\fMS\oW{n}}
corresponds by the \PSa analogue of Lemma \ref{lfreema} (the ordinary Yoneda embedding)
to an element \w{\tilde{\jmath}\sb{n}\in V\sb{n}\lin{\oW{n}}} \wwh that is. a
 homotopy class \w[.]{[\of{n}{n}]\in\pi\sb{0}\bX\sb{n}=\pi\sb{0}\fV\sb{n}\lin{\oW{n}}}
 Since $\bY$ is fibrant in $\C$, \w{\fX=\fMS\bY} is fibrant in \w[,]{\Sa\sp{\TuA}} so
 \w{\co{\fX}} is Reedy fibrant. But \w{\fVd\to\co{\fX}} is a Reedy fibration, so \w{\fVd}
  is Reedy fibrant, and therefore \w{\Xd} is, too. Thus by \cite[Lemma 2.7]{StoV}
 the inclusion of the (horizontal) Moore object
\w{C\sb{n}\Xd:=C\sb{n}\sp{h}\Xd\hra\bX\sb{n}} induces an isomorphism
\w{\pi\sb{0}C\sb{n}\Xd\to C\sb{n}\pi\sb{0}\Xd=C\sb{n}\Vd\lin{\oW{n}}}
 (see also \cite[Lemma 1.30]{BJTurHI}).

The functor \w{\fMS} of \S \ref{drdma} takes any pointed limit in $\C$ to the
corresponding limit of $\uS$-presheaves, so
\w[,]{C\sb{n-1}\fMS\WW{n-1}=\fMS C\sb{n-1}\WW{n-1}} and thus the attaching map
\w{d\sb{0}=\odz{W}{n}:\oW{n}\to C\sb{n-1}\WW{n-1}} corresponds
 under Lemma \ref{lfreema} to an element
\w[.]{\gamma\in C\sb{n-1}\fMS\WW{n-1}\lin{\oW{n}}}
Moreover, the given map of $\uS$-presheaves \w{\fff\bp{n-1}} induces
 \w[,]{C\sb{n-1}\fff\bp{n-1}:C\sb{n-1}\fMS\WW{n-1}\to C\sb{n-1}\fVd}
 which takes $\gamma$  to an element
\w[.]{\psi\sb{n}:=\fff\bp{n-1}(\gamma)\in(C\sp{h}\sb{n-1}\Xd)\sb{0}}

 Since \w{\Xd} is Reedy fibrant, the matching structure map
 \w{\delta\sb{n}:\bX\sb{n}\to M\sb{n}\Xd} is a fibration (cf.\ \cite[\S 16.3]{PHirM}),
 and we have an inclusion \w[,]{\iota:C\sb{n-1}\Xd\hra M\sb{n}\Xd}
 given by \w[.]{x\mapsto(x,x,0,\dotsc,0)} Because \w{\Wd} realizes the  CW resolution
 \w{\Vd\to\Lambda} of \PAa[s] and \w{j\sb{n}:\oV{n}\hra V\sb{n}}
factors through \w[,]{C\sb{n}\Vd} we have
\w[.]{(\delta\sb{n})\sb{\ast}[\of{n}{n}]=[\iota\circ\psi\sb{n}]}
 We may therefore change \w{\of{n}{n}} within its homotopy class so that
 \w{\delta\sb{n}(\of{n}{n})=\iota\circ\psi} on the nose.

Lemma \ref{lfreema}, together with \wref{eqrsocp} (and our assumption that
\w{\Wn{k}{n-1}} and \w{\CsW{n-k-1}{n}} are in \w[)]{\TsS} imply that
\w{\fMS\tWn{n}{n}} is the coproduct of \w{\fMS\Wn{n}{n-1}} and \w[.]{\fMS\oW{n}}
Therefore, this choice of \w{\of{n}{n}} defines a map of $\uS$-presheaves
\w{\fff\bp{n}:\fMS\tWn{n}{n}\to\fV\sb{n}} (extending \w[).]{\fff\bp{n-1}}
Since \w{F\sb{n-1}\rest{\oW{n}}=d\sp{h}\sb{0}(\gamma)} (in the notation
of \S \ref{sseqreal}(iii)), we have \w[\vsm.]{d\sp{h}\sb{0}\of{n}{n}=d\sp{h}\sb{1}\of{n}{n}}

 \noindent\textbf{Step C.}\ In the $k$-th stage of the inner (downward) induction,
 with \w[,]{k<n} we assume that for each for \w{k<j\leq n} we have chosen a map
 of \wSma[s] \w[,]{\of{n}{j}:\fMS\CsW{n-j-1}{n}\to\fV\sb{j}} represented by an element
 \w{\psi\sb{j}\in\fV\sb{j}\lin{\CsW{n-j-1}{n}}\sb{0}} with \w{d\sb{i}\sp{h}\psi\sb{j}=0}
 for \w[.]{2\leq i\leq j} If \w{\iota\sb{n-j-1}:\SW{n-j-1}{n}\hra\CsW{n-j-1}{n}}
is the inclusion, then \w{\varp{j}:=\iota\sb{n-j-1}\sp{\ast}\psi\sb{j}} lies in
\w[,]{C\sp{h}\sb{j-1}\fVd\lin{\SW{n-j-1}{n}}\sb{0}} and by induction it represents
\begin{myeq}\label{eqfj}
  \fMS\CsW{n-j-1}{n}~\xra{(F\sb{j})\sb{\ast}}~
  C\sb{j-1}\sp{h}\fMS\WW{n-1}~\xra{C\sb{j-1}\fff\bp{n-1}}~
 C\sb{j-1}\sp{h}\fVd
\end{myeq}
\noindent (in the notation of \wref[).]{eqakfk}
If \w{q\sb{n-j-2}:\CsW{n-j-2}{n}\to\SW{n-j-1}{n}} is the quotient map,
this implies that \w{q\sb{n-j-2}\sp{\ast}\varp{j}} represents
\begin{myeq}\label{eqaj}
 \fMS\Sigma\sp{n-j-1}\oW{n}~\xra{(a\sb{j})\sb{\ast}}~ Z\sb{j-1}\sp{h}\fMS\WW{n-1}~
 \xra{Z\sb{j-1}\fff\bp{n-1}}~Z\sb{j-1}\sp{h}\fVd~,
\end{myeq}
\noindent (again using the notation of \wref[),]{eqakfk} so
\w{q\sb{n-j-2}\sp{\ast}\varp{j}} is
in \w[,]{Z\sp{h}\sb{j-1}\fVd\lin{\CsW{n-j-2}{n}}\sb{0}}

Similarly, \w{d\sp{h}\sb{0}\varp{j}}
actually  lies in \w[,]{Z\sp{h}\sb{j-2}\fVd\lin{\CsW{n-j-2}{n}}\sb{0}} and represents
\begin{myeq}\label{eqajm}
 \fMS\Sigma\sp{n-j}\oW{n}~\xra{(a\sb{j-1})\sb{\ast}}~ Z\sb{j-2}\sp{h}\fMS\WW{n-1}~
 \xra{Z\sb{j-2}\fff\bp{n-1}}~Z\sb{j-2}\sp{h}\fVd~.
\end{myeq}

The nullhomotopy \w{F\sb{k}} for \w{v\sb{k-1}\circ a\sb{k}} (cf.\ \wref[).]{eqakfk}
is represented by \w[,]{\varp{k}\in C\sp{h}\sb{k-1}\fVd\lin{\SW{n-k-1}{n}}\sb{0}}
and as in Step B we use the embedding of \w{C\sp{h}\sb{k-1}\fVd\lin{\SW{n-k-1}{n}}}
in \w{M\sb{k}\fVd\lin{\SW{n-k-1}{n}}} and the facts that
\w{\delta\sb{k}:\fV\sb{k}\lin{\SW{n-k-1}{n}}\to M\sb{k}\fVd\lin{\SW{n-k-1}{n}}}
is a fibration, and that \w{\varp{k}} lifts up to homotopy to
\w{\fV\sb{j}\lin{\CsW{n-k-1}{n}}} (since \w{\CsW{n-k-1}{n}}
is contractible) to obtain an element
\w{\psi\sb{k}} in \w{\fV\sb{j}\lin{\CsW{n-k-1}{n}}\sb{0}}
(with \w{d\sb{i}\sp{h}\psi\sb{k}=0}
for \w[).]{2\leq i} such that \w[\vsm.]{\varp{k}:=d\sp{h}\sb{0}\psi\sb{k}}

\noindent\textbf{Step D.}\ The three conditions \eqref{eqfj}-\eqref{eqaj}-\eqref{eqajm}  on
 \w{\varp{j}:=d\sp{h}\sb{0}\vpsi{j}} \wb{0\leq j\leq n} are all that is needed
in order for the elements \w{\vpsi{j}} to fit together to define a map
of restricted simplicial $\uS$-presheaves \w{\tff{n-1}:\fMS\tWd{n}\to i\sp{\ast}\fVd}
extending \w{i\sp{\ast}\fff\bp{n-1}}
(in the notation of \S \ref{snac}), and so, using \wref[,]{eqrsopo} an induced map of
simplicial $\uS$-presheaves \w[,]{\hff{n-1}:\fMS\vWu{n}\to \fVd} which is a levelwise weak
equivalence through dimension $n$.

Recall from \cite[2.C]{BJTurHI} that \w{\WW{n}} is constructed by the
following factorizations in the Reedy model category structure on
\w{\Sa\sp{\TsS\op\times\Dop}}  (see \S \ref{smcsps}):
\mydiagram[\label{eqwnmon}]{
  \WW{n-1} \ar@{->}[rr] \ar@{^{(}->}[d]\sb{\irn{n-1}} &&
  \vWu{n} \ar@{^{(}->}[d]\sp{\simeq} \ar[rrd] && \\
  \WW{n} \ar@{->>}[rr]\sp{\simeq}\sb{h}  && \Wp{n} \ar@{->>}[rr] && \ast
}
\noindent where \w{\hra} indicates a cofibration and \w{\epic} a fibration,
with the top horizontal map a levelwise weak equivalence in simplicial dimensions
\www[,]{\leq n-1} so the same is true of the left vertical map.

Applying \w{\fVd\lin{-}} to \wref{eqwnmon} yields a diagram of
bisimplicial spaces, and taking diagonals, a similar diagram in
\w[.]{\Sa\sp{\Dop}} Since by our initial assumption all objects of \wref[,]{eqwnmon}
in each simplicial dimension, are in \w[,]{\TsS}, by
Lemma \ref{lfreema} we obtain an analogous diagram of mapping
spaces of $\uS$-presheaves into \w[.]{\fVd} The sequence of
elements in the simplicial set \w{\diag\fVd\lin{\WW{n-1}}\sb{0}}
in the upper left corner corresponding to
\w{\fff\bp{n-1}:\fMS\WW{n-1}\to\fVd} map by construction to the
sequence in \w{\diag\fVd\lin{\vWu{n}}\sb{0}} corresponding to
\w[,]{\hff{n-1}} mapping forward to a sequence $\beta$
corresponding to \w[.]{\ppp\fff\bp{n}:\fMS\Wp{n}\to\fVd} Since the
map $h$ in \wref{eqwnmon} is a trivial fibration, and these are
preserved by evaluation of \w{\fVd} and diagonals, we see that the
induced map of simplicial sets
\w{h\sb{\ast}:\diag\fVd\lin{\WW{n}}\sb{0}\to\diag\fVd\lin{\Wp{n}}\sb{0}}
is a trivial fibration. We can therefore lift $\beta$ to a
sequence representing the required map
\w[,]{\fff\bp{n}:\fMS\WW{n}\to\fVd} completing the outer induction
step.
\end{proof}

\begin{remark}\label{rtlwe}
The same result holds if we replace $\uS$-presheaves by $r$-truncated $\uS$-presheaves,
since (as noted in \S \ref{struncma}), Lemma \ref{lfreema} still holds, and
\w{\Wd:=\Po{r}\fMS\Wd} is free in each simplicial dimension.
\end{remark}

\begin{summary}\label{srldf}
Assume given a homotopy functor \w[,]{\bT:\G\to\M} inducing
\w{\fT:=\bT\circ N:\sMst\to\M} as in Corollary \ref{crssma}.
Let \w{\D:=\Sa\sp{\TsSt\op\times\Dop}} and \w[,]{\E:=\M\sp{\Dop}}
with the resolution model category structure on $\D$ determined by $\uS$ for $\G$
as in \S \ref{smcsps}, with respect to the structure
of \S \ref{smcma} for \w{\Sa\sp{\TsSt\op}} (with \ww{E\sp{s}}-weak
equivalences on $\E$).

In the notation of \S \ref{srdf}, let $\C$ denote the category of
simplicial \sSma[s] in $\D$ associated to sequential realizations as in
\S \ref{sseqreal}, let \w{\cW'} be the Reedy weak equivalences in $\D$, and let
\w{\D\sb{\C}} be the full subcategory \w{\ho\sMst} of objects in
\w{\ho(\Sa\sp{\TsSt\op\times\Dop})} weakly equivalent
to a constant simplicial object on \w[.]{\sMst}
The \emph{relative left derived functor} \w{\bLr\fT:\ho\sMst\to\ho\E}
is then defined on an \Stma $\fX$ (more formally, on \w[)]{\co{\fX}} by
\begin{enumerate}
\renewcommand{\labelenumi}{(\alph{enumi})}
\item Choosing a simplicial resolution \w{\bet:\fVd\to\fX} in
\w[;]{\Sa\sp{\TsS\op\times\Dop}}
\item Choosing a CW basis \w{\{\oV{n}\}\sb{n=0}\sp{\infty}} for the \PSa[-resolution]
\w[,]{\Vd:=\pi\sb{0}\fVd\to\pi\sb{0}\fX} a sequential realization $\cW$ of \w{\Vd}
for \w[,]{\bY:=N\fX} with an \ww{E\sp{1}}-weak equivalence \w[;]{\fMst\Wd\to\fVd}
\item Defining \w{(\bLr\fT)\fX} to be the simplicial object \w{\fT\fMst\Wd} in
\w{\ho\D\sp{\Dop}} (uniquely determined up to \ww{E\sp{1}}-weak equivalence).
\end{enumerate}
\end{summary}

\supsect{\protect{\ref{chotdf}}.B}{Relative derived functors of dual mapping algebras}

For a given commutative ring $R$, let \w{\TuR=\sMRla} be the full subsimplicial
category of \w{\C=\SR} consisting of all simplicial $R$-modules of cardinality
$<\lambda$, as in \S \ref{detr}, and \w{\fX=\fMR\bY} for some \w{\bY\in\C}
(the cardinal $\lambda$ we choose may depend on $\bY$).
Essentially, we may dualize the results of \S \ref{chotdf}.A to this situation.
Note that because \w{\bY\mapsto H\sp{\ast}(\bY;R)} is contravariant
the category \w{\PAAlg} resembles \w{\PSAlg} in being a category of graded
universal algebras,  so the resolutions we need for the \PAa
\w{\Lambda=H\sp{\ast}(\bY;R)} will be simplicial, rather than cosimplicial,
and we can use the notion of a CW resolution \w{\Vd\to\Lambda} as in
\S \ref{scompres}.  However, only when $R$ is a field do we know that any
free simplicial resolution in \w{\PAAlg\sp{\Dop}} has a CW basis
\w{(\oV{n})\sb{n=0}\sp{\infty}} of free \PAa[s] (see \cite[Proposition 3.12]{BlaS}).
For the cosimplicial resolutions of spaces, we need to dualize \S \ref{scompres}
as follows:

\begin{defn}\label{dmcoch}
If $\C$ is cocomplete, the $n$-th \emph{Moore cochain} object of a cosimplicial
object \w{\Gu\in\C\sp{\Delta}} is
\w[,]{C\sp{n}\Gu:=
\Cok(\coprod\sb{i=1}\sp{n-1}\,G\sp{n}\xra{\bot\sb{i}\,d\sp{i}}G\sp{n})}
with differential \w{\delta\sp{n-1}:C\sp{n-1}\Gu\to C\sp{n}\Gu} induced by
\w[,]{d\sp{0}\sb{n-1}} and structure map \w[.]{v\sp{n}:G\sp{n}\to C\sp{n}\Gu}
We denote the cofiber of \w{\delta\sp{n-1}} by \w[,]{Z\sp{n}\Gu}
with structure map \w[,]{w\sp{n}:C\sp{n}\Gu\epic Z\sp{n}\Gu} and note that
\w{\delta\sp{n-1}} factors as \w[.]{\udz{n-1}\circ w\sp{n-1}}
\end{defn}

\begin{mysubsection}{Dual sequential realizations}\label{sdseqreal}
Let $R$ be a commutative ring and $\lambda$ a limit cardinal, with \w{\ThuA:=\pi\sb{0}\TuR}
  for \w[.]{\TuR=\sMRla} Assume given  an $R$-good space \w{\bY\in\Sa} and a
  CW resolution \w{\Vd} of the \PAa \w[,]{\Lambda=\piA\bY} with CW basis
  \w[,]{\{\oV{n}\}\sb{n=0}\sp{\infty}} such that for each \w[,]{n\geq 0}
  \w{\oV{n}\cong\piA\uW{n}} for some \w[.]{\uW{n}\in\TuR}

  We define a (dual) \emph{sequential realization of \w{\Vd} for $\bY$} to be a sequence
  $\cW$ of maps
\begin{myeq}\label{eqdtower}
\dotsc~\WWu{n+1}~\xra{\prn{n+1}}~\WWu{n}~\xra{\prn{n}}~\WWu{n-1}~\dotsc~
\WWu{1}~\xra{\prn{1}}~\WWu{0}
\end{myeq}
\noindent between Reedy fibrant and cofibrant objects in
\w[,]{\Sa\sp{\Delta}} such that for each \w[:]{n\geq 0}

\begin{enumerate}
\renewcommand{\labelenumi}{(\roman{enumi})~}
\item There is an $n$-skeletal restricted cosimplicial object \w{\tWu{n}} with
  \w{\tWun{k}{n}=\Wun{k}{n-1}\times\PoW{n-k-1}{n}} for \w[,]{0\leq k\leq n}
 where as before by convention \w[.]{\ooW{0}{n}=\PoW{-1}{n}=\uW{n}}
\item The coface map \w{d\sp{0}:C\sp{k}\to\tWun{k+1}{n}} into the factor
  \w{\PoW{n-k-2}{n}} is the map \w{F\sp{k}} in the commuting diagram
\mydiagram[\label{eqdakfk}]{
  Z\sp{k-1}\WWu{n-1} \ar[rr]\sp{\udz{k-1}} \ar[d]\sp{a\sp{k-1}} &&
  C\sp{k}\WWu{n-1} \ar@{->>}[rr]\sp{w\sp{k}} \ar[d]\sp{F\sp{k}} &&
Z\sp{k}\WWu{n-1} \ar[d]^{a\sp{k}} \\
  \ooW{n-k-1}{n} \ar@{^{(}->}[rr]^{j\sp{n-k-1}} &&   \PoW{n-k-2}{n} \ar@{->>}[rr]^{p\sp{n-k-2}} &&
  \ooW{n-k-2}{n}
}
\noindent (in the notation of \S \ref{dmcoch}).
The first coface map \w{d\sp{1}} into \w{\PoW{n-k-2}{n}} is the composite
of  the projection onto \w{\PoW{n-k-1}{n}} with \w[,]{j\sp{n-k-1}\circ p\sp{n-k-1}}
and \w{d\sp{i}} into the factor  \w{\PoW{n-k-2}{n}} is zero for \w[.]{i>1}

We start with a realization of the $n$-th
attaching map \w{\odz{V}{n}:\oV{n}\to C\sb{n-1}\Vd} for the given CW resolution
as our choice for
\w[.]{F\sp{n-1}:C\sp{n-1}\WWu{n-1}\to\uW{n}}
\item Let  \w{\vWun{n}} be the pullback of
\w[,]{\WWu{n-1}~\leftarrow\F i\sp{\ast}\WWu{n-1}~\to~\F\tWd{n}}
where \w{\F:\C\sp{\Delta}\to \C\sp{\Delta}} is the right adjoint of
the forgetful functor \w{i\sp{\ast}:\C\sp{\Delta}\to \C\sp{\Dz}} (see \S \ref{snac}),
with \w{\WWu{n}} a Reedy fibrant and cofibrant replacement for \w[.]{\vWun{n}}
\end{enumerate}

Again, \w{\Wu:=\lim\sb{n}\WWu{n}} is a cosimplicial resolution of $\bY$
in the resolution model category \w[,]{\C\sp{\Delta}}  and in fact, the sequential
realization $\cW$ can be constructed starting from any \Rma $\fX$.
See \cite[\S 2 \& Appendix A]{BSenH} for further details.
\end{mysubsection}

The proof of Theorem \ref{ttlwe} can be dualized to yield:

\begin{thm}\label{tdtlwe}
  Given a commutative ring $R$ with \w{\TuA=\sMRla} and an $R$-good space $\bY$, let
  \w{\bet:\fVd\to\fX=\fMR\bY} be a simplicial resolution in \w{\Sa\sp{\TuR\times\Dop}}
with a CW basis \w{\{\oV{n}\}\sb{n=0}\sp{\infty}} for the \PAa[-resolution]
\w[.]{\Vd:=\pi\sb{0}\fVd\to\Lambda=\piA\bY} Then for any sequential realization
$\cW$ of \w{\Vd} for $\bY$, there is a Reedy weak equivalence of simplicial
\wAma[s] \w[.]{\fff:\fWd:=\fMA\Wu\to\fVd}
\end{thm}

The dual of Remark \ref{rtlwe}, for the $n$-truncated case, also holds.

\begin{summary}\label{ddthodf}
Given a functor \w{\fT:\dsMsr\to\D} as in Corollary \ref{cdrssma},
the \emph{relative right derived functor}
\w{\bRr\fT:\ho\dsMsr\to\ho(\D\sp{\Delta})}
applied to \w{\fX:=\fMRst\bY} for $R$-good \w[,]{\bY\in\Sa} is obtained by
\begin{enumerate}
\renewcommand{\labelenumi}{(\alph{enumi})}
\item Choosing a simplicial resolution \w{\bet:\fVd\to\fX} in the model
  category \w[;]{\Sa\sp{\TuA\times\Dop}}
\item Assuming the \PAa[-resolution] \w{\Vd:=\pi\sb{0}\fVd\to\piA\bY}
has a CW basis \w{\{\oV{n}\}\sb{n=0}\sp{\infty}} (e.g., if $R$ is a field),
choosing a sequential realization $\cW$ of \w[;]{\Vd}
\item Defining \w{(\bRr\fT)\fX} to be the cosimplicial object
\w{\fT\fMRst\Wu} in \w[.]{\D\sp{\Delta}}
\end{enumerate}
\end{summary}

\begin{example}\label{egdthodf}
For \w{\bZ\in\Sa} and \w{\bT:=\mapa(\bZ,-)} as in \S \ref{egdrssma},
if \w{\fZ=\fMRst\bZ} and \w{\fX=\fMRst\bY} for some $R$-good space $\bY$,
and \w{\fVd=\fMRst\Wu} for some cosimplicial resolution \w[,]{\bY\to\Wu}
then \w{(\bRr\fT)\fX:=\fT\fVd} is the cosimplicial space \w{\fZ\lin{\Wu}} (up to
\ww{E\sp{2}}-equivalence).
\end{example}

%
%
\sect{Truncating higher order derived functors}
\label{cthodf}

So far we have shown only that the usual total derived functor
\w{\bLL\bT} of a continuous functor \w{\bT:\C\to\D} can be
interpreted (under suitable assumptions) as derived functors of
the corresponding \ma[s.] Although there are many technicalities
involved, the result is hardly surprising, since, under these
assumptions, \ma[s] carry the same homotopy information as objects
in $\C$  (Theorems \ref{trssma} and \ref{tdrssma}).

The point is that \ma[s] are the right framework for \emph{truncating} the
homotopy information (using Postnikov sections), while still retaining enough
to compute the required term in the homotopy spectral sequences
for \w{\bT\Wd} or \w[.]{\bT\Wu}

\supsect{\protect{\ref{cthodf}}.A}{Truncating derived functors of \ma[s]}

Not every homotopy functor $\bT$ (and the corresponding $\fT$) will behave
as we want with respect to such truncation. We therefore require the
following:

\begin{defn}\label{dtrunss}
For any \w[,]{2\leq r\leq \infty} let \w{\E\sp{r}} denote the category of
$r$-truncated homological spectral sequences
\w[,]{\{E\sp{k}\sb{\ast\ast}\}\sb{k=1}\sp{r}} equipped with a differential
\w[,]{d\sp{r}:E\sp{r}\sb{t,i}\to E\sp{r}\sb{t-r-1,i+r}}
which need not satisfy \w[.]{d\sp{r}\circ d\sp{r}=0} A map in \w{\E\sp{r}}
is called a \emph{weak equivalence} if it induces an isomorphism in
\w{E\sp{2}\sb{\ast\ast}} (and thus also for \w[).]{r\leq k>2} This defines the
corresponding localized category \w[.]{\ho\E\sp{r}} We have truncation
functors \w{\Po{r}:\E\sp{n}\to\E\sp{r}} for each \w[.]{r\leq n\leq \infty}  Note that
the homotopy spectral sequence of a simplicial space defines a
homotopy functor \w{\eS\sp{\infty}:\G\sp{\Dop}\to\E\sp{\infty}} (with respect to
\ww{E\sp{2}}-equivalences in the source and target), and write
\w[.]{\eS\sp{r}:=\Po{r}\circ\eS\sp{\infty}}
\end{defn}

\begin{defn}\label{dlevelf}
Any homotopy functor \w[,]{\bT:\G\to\G} and the corresponding
\w[,]{\fT:\sMst\to\G} induce a functor
\w{\eS\sp{r}\circ\bLr\fT:\ho\sMst\to\ho\E\sp{r}} (see \S \ref{srldf}) for each
\w[.]{r\geq 2} We say that $\bT$ (and $\fT$) are \emph{level} if for
every \w[,]{r\geq 2} this functor \w{\eS\sp{r}\circ\bLr\fT}
factors through a functor \w[.]{\bLr\fTn{r-2}:\ho\sMstn{r-2}\to\ho\E\sp{r}}

Here \w{\ho\sMstn{n}} is the subcategory of \w{\ho(\Sn{n}\sp{\TsSt\op\times\Dop})}
weakly equivalent to \w[,]{\co{\fX}} for $\fX$ in the subcategory
\w{\sMstn{n}} of $n$-truncated \Stma[s] (cf.\ \S \ref{struncma}).
\end{defn}

In order to identify which homotopy functors are level, we shall need the
following notion introduced in \cite[\S 1]{BBlaS} (see also \cite{BDGoeR}):

\begin{defn}\label{dpstem}
Let $\C$ be \w[,]{\Tz} \w[,]{\Sa} or $\G$: for any \w[,]{n\geq 0} an \emph{$n$-stem}
in $\C$ is a tower:
\begin{myeq}[\label{eqponstem}]
\cQ~:=~\left(\dotsc~\to~Q\sb{k+1}~\xra{q\sb{k+1}}~Q\sb{k}~\xra{q\sb{k}}~
Q\sb{k-1}~\dotsc~Q\sb{1}\right)
\end{myeq}
\noindent in \w[,]{\C\sp{(\NN,\leq)}} in which
\w{\pi\sb{i}(Q\sb{k})=0} for \w{i<k} or \w[,]{i>n+k} and
\w{\pi\sb{i}q\sb{k}} is an isomorphism for \w[.]{k\leq i<n+k}
Here \w{(\NN,\leq)} is the usual linearly ordered category of the natural
numbers. The object \w{Q\sb{k}\in\C} is called the $k$-th \emph{$n$-window}
of $\cQ$.

We denote by \w{\Stem[n]} the full subcategory of $n$-stems
in the functor category \w[,] {\C\sp{(\NN,\leq)}} with the model
category structure on the latter as in \cite[11.6]{PHirM}. The
\emph{Postnikov $n$-stem} functor \w{\PPo{n}:\C\to\Stem[n]} is given by
\w[.]{\PPo{n}X:=\{\Po{n+k+1}X\lra{k}\}\sb{k=1}\sp{\infty}}

To avoid the need to distinguish the cases \w{\C=\Tz} or $\G$, we everywhere
use the \ww{\Top}-indexing for spheres, homotopy groups, Postnikov systems,
and connected covers (as in \S \ref{sstovcons}).
\end{defn}

By \cite[Theorem 4.13 \& Corollary 4.16]{BBlaS} we have:

\begin{thm}\label{tsss}
For each \w{r\geq 2} there is a functor \w{\hS{r}:\Stem[r-1]\sp{\Dop}\to\E\sp{r}}
which associates to any simplicial \wwb{r-1}stem \w{\Qd} an
$r$-truncated spectral sequence. Moreover,
\w{\hS{r}\circ\PPo{r-1}:\C\sp{\Dop}\to\E\sp{r}} is naturally equivalent to
\w[,]{\eS\sp{r}} so when \w{\Qd=\PPo{r-1}\Xd} this is the truncation of the usual homotopy
spectral sequence for \w[.]{\Xd} In this case we have
\w[,]{d\sp{r}\circ d\sp{r}=0} so in fact the spectral sequence is determined through
\w{E\sp{r+1}\sb{\ast\ast}} (though without \w[).]{d\sb{r+1}}
\end{thm}

\begin{cor}\label{csss}
A functor \w{\fT:\sMst\to\G} associated to a homotopy functor \w{\bT:\G\to\G}
is level if for each \w[,]{r\geq 1} the relative derived
functor \w{\eS\sp{r}\circ\bLr\fT:\ho\sMst\to\E\sp{r}}
factors as \w{\hS{r}\circ\bLr\fTn{r-1}} for some functor
\w[.]{\bLr\fTn{r-1}:\ho\sMstn{r-1}\to\ho(\Stem[r-1]\sp{\Dop})}
\end{cor}

In order for Corollary \ref{csss} to be of any use, we must identify level
homotopy functors $\bT$ for which the homotopy spectral sequence of \w{\bT\Xd}
is of interest. We first note:

\begin{lemma}\label{ltrunstem}
  For $\uS$ as in \S \ref{egmapalg}, any $n$-truncated \wSma \w{\fX\in\sMstn{n}} is
  functorially realizable  by an $n$-stem
\w[.]{\cQ=\{Q\sb{k}\}\sb{k=1}\sp{\infty}} Moreover, if \w{\fX=\Po{n}\fMS\bY} for some
\w[,]{\bY\in\G} then $\cQ$ is naturally weakly equivalent to the Postnikov
$n$-stem \w[.]{\PPo{n}\bY}
\end{lemma}

\begin{proof}
This result appears in \cite[\S 10.5]{BBlaC} for \Stma[s,] but in fact we need
only observe that for \w[,]{n\geq 1} the action of \w{\Po{n}\TsS} on $\fX$
includes \emph{inter alia} an \ww{A\sb{n}}-structure on
\w[,]{X\sb{k}:=\fX\lin{\bS{k}}} so allowing \w{\Po{n}X\sb{k}}
to be delooped to produce the window \w{Q\sb{k}} by \cite[Corollary 11.12]{StasH}.
The weak equivalences \wref[,]{eqtruncloop} together with \cite[Theorem 12.7]{MayG},
yield the structure maps for the $n$-stem $\cQ$
\end{proof}

The simplest example is from \cite{BlaH}, where it is used to construct a spectral
sequence for computing \w{H\sb{\ast}\bY} from the \Pa \w[:]{\pis\bY}

\begin{prop}\label{phure}
The abelianization functor \w{\Ab:\G\to\G} is level.
\end{prop}

\begin{proof}
Let \w{\cQ=\{Q\sb{k}\}\sb{k=1}\sp{\infty}} denote the Postnikov $n$-stem of a space $\bX$, and
\w{\R=\{R\sb{k}\}\sb{k=1}\sp{\infty}} that of
\w[.]{\Ab\bX} Note that for each \w[,]{k\geq 0} the covering map
\w{\rho:\bX\lra{k}\to\bX} induces a map \w[,]{\rho\sb{\ast}:\Ab(\bX\lra{k})\to\Ab\bX}
which factors through \w{(\Ab\bX)\lra{k}} by cellularity (uniquely, if we choose a
\wwb{k+1}reduced model for connected covers \wh which is an inclusion of a sub-simplicial group,
in $\G$). Furthermore, by the Hurewicz Theorem,
for each \w{m\geq 0} the structure map \w{p\sb{m}:\bX\to\Po{m}\bX}
induces an isomorphism \w{H\sb{i}\bX\to H\sb{i}\Po{m}\bX} for \w[,]{i\leq m} and
an epimorphism \w[,]{H\sb{m+1}\bX\to H\sb{m+1}\Po{m}\bX}
so the natural map \w{\Po{m}(\Ab\bX)\to\Po{m}\Ab(\Po{m}\bX)} is a weak equivalence.
Thus we have a natural weak
equivalence \w{\Po{n+k+1}(\Ab Q\sb{k})\simeq R\sb{k}} for each \w[.]{k\geq 0}

Thus a given a simplicial resolution \w{\fVd\to\Po{n}\fX=\Po{n}\fMS\bY}
of $n$-truncated $\uS$-presheaves in the model category \w[,]{\Sn{n}\sp{\TsS\op}}
by Lemma \ref{ltrunstem}, we obtain a simplicial $n$-stem \w[,]{\cQd} which
yields in turn the required simplicial $n$-stem \w[.]{\Rd:=\PPo{n}(\Ab\cQd)}
\end{proof}

Here are two additional examples from \cite{StoV}. The first is used to construct
a spectral sequence for computing \w{\pis\Sigma\bY} from \w[:]{\pis\bY}

\begin{prop}\label{psusp}
The suspension functor \w{\Sigma:\G\to\G} is level.
\end{prop}

\begin{proof}
For each \w[,]{n\geq 1} any $n$-truncated \wSma has a corresponding
$n$-stem $\cQ$ by Lemma \ref{ltrunstem}, and  the \Pa \w{\Lambda:=\pi\sb{0}\fX}
determines the \Pa structure on \w{\pis Q\sb{k}} for
each \w[.]{k\geq 0}  If \w{\fX\simeq\Po{n}\fMS\bX} for some space $\bX$, then $\Lambda$
is isomorphic to \w{\pis\bX} and \w[.]{Q\sb{k}\simeq\Po{n+k+1}\bX\lra{k}}
To understand \w[,]{\bLL\fT} we need only consider the case
when $\Lambda$ is a free \Pa[.]

Now let \w[,]{\R=\{R\sb{k}\}\sb{k=1}\sp{\infty}} denote the Postnikov $n$-stem of
\w[.]{\Sigma\bX} As in the proof of Proposition \ref{phure}, the covering map
\w{\rho:\bX\lra{k}\to\bX} induces a map
\w[.]{\rho\sb{\ast}:\Sigma(\bX\lra{k})\to(\Sigma\bX)\lra{k+1}}
Taking Postnikov sections yields natural maps
\w[.]{p\sb{k}:\Po{n+k+2}(\Sigma Q\sb{k})\to R\sb{k+1}}
In particular, \w{p\sb{0}:\Po{n+2}(\Sigma Q\sb{0})\to R\sb{1}=\Po{n+1}(\Sigma\bX)}
is a weak equivalence by the Hurewicz Theorem, with \w{\Po{1}R\sb{1}\simeq\bX\sb{1}}
(a wedge of $1$-spheres, and thus aspherical).

However, for \w{k>1} there is no functorial description of \w{R\sb{k}} in terms of
$\cQ$. Thus if \w{\fT:=\Sigma\circ N:\sMst\to\G} is induced by \w{\Sigma:\G\to\G}
as in Corollary \ref{crssma}, in order to define
\w{\bLL\fTn{n}:\ho\sMstn{n}\to\Stem[n]} we must proceed as follows:

By Lemma \ref{ltrunstem} a simplicial resolution \w{\fVd\to\Po{n}\fX=\Po{n}\fMS\bY}
of $n$-truncated \wSma[s] yields a simplicial $n$-stem \w[.]{\cQd} Since
the simplicial \Pa \w{\Vd:=\pi\sb{o}\fVd} is a free resolution of
\w[,]{\Lambda:=\pi\sb{0}\fX} it has a (non-canonical) CW basis
\w{\{\oV{n}\}\sb{n=0}\sp{\infty}} for it, which in turn has a sequential realization
$\cW$ (see \S \ref{sseqreal}). By Remark \ref{rtlwe},
there is a Reedy weak equivalence of simplicial $n$-truncated \wSma[s]
\w[,]{\fff:\Po{n}\fWd\to\fVd}  where \w{\fWd} is realizable as
\w[.]{\fMS\Wd} We can realize \w{\Po{n}\fWd} by the simplicial $n$-stem
\w[,]{\hQd\simeq\PPo{n}\Wd} and let \w{\Sigma\hQd} denote the simplicial $n$-stem
obtained by applying $\Sigma$ to each window of \w{\hQd} (and taking appropriate
Postnikov sections). If \w{\hRd} denotes the simplicial Postnikov $n$-stem
\w[,]{\PPo{n}\Sigma\Wd} we have a map of simplicial $n$-stems
\w[,]{\widehat{p}:\Sigma\hQd\to\hRd} as explained above.

Similarly, the simplicial $n$-truncated $\uS$-presheaf \w{\fVd} yields a simplicial
$n$-stem \w[,]{\cQd} and \w{\fff:\Po{n}\fWd\to\fVd} induces a levelwise weak equivalence
of simplicial $n$-stems \w{\widehat{f}:\Sigma\hQd\to\Sigma\cQd}
(in the Reedy model structure). We may assume that each window of all the simplicial
$n$-stems described here are cofibrant in $\G$, so they are Reedy cofibrant. Thus if
we let \w{\Rd} denote the homotopy pushout of $\widehat{f}$ and $\widehat{p}$ (in
the Reedy model category of simplicial $\uS$-presheaves), we have
a Reedy weak equivalence \w{\hRd\to\Rd} (cf.\ \cite[Proposition 13.1.2]{PHirM}),
as well as a structure map of simplicial $n$-stems \w[.]{p:\cQd\to\Rd}

We define \w{(\bLL\fTn{n})\Po{n}\fX} to be the simplicial $n$-stem
\w{\Rd} To see that \w{\bLL\fTn{n}} is well-defined, replace
\w{\fVd} by  some other simplicial resolution \w{\fUd\to\Po{n}\fX}
of $n$-truncated $\uS$-presheaves, with \w{\cZ} a sequential
realization of \w{\pi\sb{0}\fUd} for $\bY$. Let \w{\Rd} and
\w{\Sd} denote the simplicial $n$-stems associated as above to
\w{\fVd} and \w{\fUd} respectively. We then have a weak
equivalence of simplicial spaces \w{g:\Wd\to\Zd} in the resolution
model category structure with respect to \w{\TsS} (since both are
cofibrant replacements for \w[),]{\co{\bY}} and this will induce a
weak equivalence \w{\fVd\to\fUd} in the resolution model structure
of \S \ref{smcsps}, and thus the same holds for the simplicial
$n$-stems \w{\Rd} and \w{\Sd} (cf.\ \cite[Theorem 1.9]{StoV}).
\end{proof}

The next example is used to construct a van Kampen spectral sequence to
compute \w{\pis(\bY\vee\bZ)} from \w{\pis\bY} and \w[:]{\pis\bZ}

\begin{prop}\label{pwedge}
The wedge bifunctor \w{\vee:\G\times\G\to\G} is level.
\end{prop}

\begin{proof}
The proof is entirely analogous to that of Proposition \ref{psusp}: given
two \Stma[s] $\fX$ and $\fY$, realizable by $\bY$ and $\bZ$, respectively, their
$n$-truncations are realizable by $n$-stems $\cQ$ and $\cS$, weakly equivalent
to the Postnikov $n$-stem \w{\PPo{n}\bY} and \w[,]{\PPo{n}\bZ} respectively.
Once again we cannot reconstruct the Postnikov $n$-stem for \w{\bY\vee\bZ} directly
from the window-wise wedge of $\cQ$ and $\cS$ (except for the bottom window), but
must have recourse to sequential realizations of the full simplicial resolutions.
\end{proof}

\begin{remark}\label{rhocolim}
Stover set up spectral sequences for arbitrary homotopy colimits in
\w{\Tz} (see \cite[Theorem 1.2]{StoV}), and one can obtain similar results for the
left derived functors appearing as the \Eut[s]{2} of these spectral sequences.
\end{remark}

\supsect{\protect{\ref{cthodf}}.B}{Truncating derived functors of
dual mapping algebras}

We may dualize Definitions \ref{dtrunss} and \ref{dlevelf}  as follows:

\begin{defn}\label{ddlevelf}
For any \w{2\leq r\leq \infty} we let \w{\E\sb{r}}
denote the category of $r$-truncated cohomological spectral sequences
\w{\{E\sb{k}\sp{\ast\ast}\}\sb{k=1}\sp{r}} (again, the last differential
need not satisfy \w[).]{d\sb{r}\circ d\sb{r}=0} A \emph{weak equivalence}
in \w{\E\sb{r}} is a map inducing an isomorphism in \w[.]{E\sb{2}\sp{\ast\ast}}
Again we have truncation functors \w[.]{\Po{r}:\E\sb{n}\to\E\sb{r}}
The homotopy spectral sequence of a cosimplicial space defines a
homotopy functor \w[,]{\eS\sb{\infty}:\Sa\sp{\Delta}\to\E\sb{\infty}} and we write
\w[.]{\eS\sb{r}:=\Po{r}\circ\eS\sb{\infty}}

If \w{\bT:\SR\to\D} is a homotopy functor preserving $R$-equivalences, we say
that $\bT$, and the corresponding \w{\fT:\dsMsr\to\D} of Corollary \ref{cdrssma},
are \emph{level} if for any \w{r\geq 2} and weakly $R$-good \dsStma
\w{\fX=\fMRst\bY} (see \S \ref{dewr}), \w{\eS\sb{r}\bRr\fT\fX}
factors up to isomorphism through the
\wwb{r-2}truncated simplicial \dsStma \w[,]{\Po{r-2}\fMRst(\bRr\fT\fX)} up to
weak equivalence in \w{\Sn{n}\sp{\TuA\times\Dop}} (see \S \ref{rmcspa}).
\end{defn}

Although the analogue of Theorem \ref{tsss} was also shown in \cite{BBlaS} to hold
for the homotopy spectral sequence of a cosimplicial space, this does not appear
to be helpful in showing that functors of \Rma[s] are level \wh mainly because there
is no simple connection between maps into Eilenberg-Mac~Lane spaces and maps out
of spheres. Thus a more direct approach is needed here.

Our main result in this connection, which may be of independent interest, is the
following reinterpretation of the results of \cite{BBSenH}:

\begin{thm}\label{tdlevelf}
For any \w{\bZ\in\Sa} and \w{R=\Fp} or $\QQ$, the unstable $R$-Adams spectral sequence
for \w{\bT:=\mapa(\bZ,-)} applied to \w{\bY\in\SR} (see \cite[\S 7.2]{BKanS})
is determined by the simplicial \Rma \w[,]{(\fMR\bRr\fT)\fMRst\bY} and $\bT$ is level.
\end{thm}

\begin{proof}
Let \w{\bY\to\Wu} be a cosimplicial resolution, which
we may assume without loss of generality to be associated to a dual
sequential realization $\cW$ as in \S \ref{chotdf}.B, by Definition \ref{ddthodf}.

We know that the homotopy spectral sequence for the cosimplicial space
\w{\Xu:=\map(\bZ,\,\Wu)} is determined in principle by the simplicial
\sAma \w[.]{\fWd:=\fMRst\Wu} Following the description in \cite{BBSenH} (and compare
\cite{BousH}) we now explain how this can be made explicit:

By \cite[Proposition 4.18]{BBSenH} the unstable Adams spectral sequence for $\bY$
as above agrees from the \Elt{2} on with that associated to the fibration sequences
\begin{myeq}[\label{eqfibseqsr}]
\Omega\sp{n}\uW{n}~\to~\Tot\sb{n}\dvWu{n}~\to~\Tot\sb{n-1}\WWu{n-1}~,
\end{myeq}
\noindent in the notation of \S \ref{sdseqreal}, so the same is true of the
homotopy spectral sequence for \w[,]{\Xu:=\map(\bZ,\,\Wu)}
if we apply \w{\mapa(\bZ,-)} before taking \w[.]{\Tot}

An element \w{\gamma\in E\sb{1}\sp{n,k+n}} is thus represented by a map
\w[,]{\Sigma\sp{k}\bZ\to\Tot\sb{n}\sDu{n}} where \w{\sDu{n}} is the fiber
of the Reedy fibration \w{\dvWu{n}\to\WWu{n-1}} and
\w{\Tot\sb{n}\sDu{n}\simeq \Omega\sp{n}\uW{n}} (see \cite[Proposition 4.12]{BBSenH}).
This is represented in turn by a map of cosimplicial spaces
\w[:]{\Gu:\bDeu\ltimes\Sigma\sp{k}\bZ\to\WWu{n}} (see \S \ref{sdseqreal}(iii)) \wh
that is, a sequence of maps
\w{G\sp{j}\lp{n}:\bDel\sp{j}\ltimes\Sigma\sp{k}\bZ\to\Wun{j}{n}} (where we may assume
\w{G\sp{j}\lp{n}=0} for \w{j<n} by \cite[(3.6)]{BBSenH}).

By \cite[Theorem 5.9]{BBSenH}, for each \w{r\geq 2} and \w[,]{N:=n+r-1} the differential
\w{d\sb{r}:E\sb{r}\sp{n,k+n}\to E\sb{r}\sp{N+1,k+N}}
is defined on \w{\lra{\gamma}} by the value
\w{\phi:\Sigma\sp{k}\bZ\to\Omega\sp{N}\uW{N+1}} of a certain $r$-th order $R$-cohomology
operation. This operation is defined when the associated sequence of lower order
operations vanish, so that there exists a chosen lift of \w{\Gu} to
\w[.]{G\lp{N}:\bDeu\ltimes\Sigma\sp{k}\bZ\to\WWu{N}}

The map $\phi$ is obtained by patching together the composite of the maps \w{G\sp{i}\lp{N}}
with the given maps \w{F\sp{j}\lp{N+1}:\Wun{j}{N}\to P\Omega\sp{N-j-1}\uW{N+1}} of
\wref[.]{eqdakfk} yielding a map from the boundary of a certain
\wwb{N+1}dimensional polyhedron \w[,]{\PP\sp{N+1}\sb{r}} described in
\cite[\S 4.3]{BBSenH} to \w[.]{\mapa(\Sigma\sp{k}\bZ,\,\uW{N+1})} This is adjoint to
a map \w[,]{\widetilde{\phi}:\Sigma\sp{k}\bZ\to\Omega\sp{N}\uW{N+1}}
and by \cite[Theorem 5.10]{BBSenH}, the class
$$
[\widetilde{\phi}]\in~[\Sigma\sp{k}\bZ,\,\Omega\sp{N}\uW{N+1}]~\cong~
[\Sigma\sp{k-1}\bZ,\,\Omega\sp{N+1}\uW{N+1}]~\cong~
E\sb{1}\sp{N+1,k+N}
$$
\noindent (using the usual $\Sigma$-$\Omega$ adjunction on the left) represents
\w[.]{d\sb{r}\lra{\gamma}\in E\sb{r}\sp{N+1,k+N}}
In particular, by \cite[Lemma 5.7]{BBSenH}, \w{[\widetilde{\phi}]} vanishes if and only
if \w{G\lp{N}} lifts to a map \w[.]{G\lp{N+1}:\bDeu\ltimes\Sigma\sp{k}\bZ\to\WWu{N+1}}

Because we assumed that each \w{\uW{N}} is in \w{\TuR} (see \S \ref{sdseqreal}),
the information used in defining this higher operation is encoded by
\w{\fWd:=\fMR\Wu} and \w{\fZ:=\fMR\bZ}
Furthermore, since \w{G\sp{j}\lp{N}=0} for \w[,]{j<n} and \w{\WWu{N}}
is \wwb{n+r-1}skeletal by \S \ref{sdseqreal}(i), from the description
above we see that we only need \w{\Po{r-1}\fZ\lin{\Omega\sp{k}\uW{N}}}
in order to calculate \w[,]{d\sb{r}} and thus \w[.]{E\sb{r+1}\sp{\ast\ast}}
Finally, by \S \ref{egdrssma}, \w{\Po{r-1}\fZ} is completely determined by
the \wwb{r-1}truncated \Rma \w[,]{\Po{r-1}\fWd} and this in turn depends only on
\w[,]{\Po{r-1}\fMRst\bY} up to \ww{E\sb{2}}-equivalence.
\end{proof}

\begin{cor}\label{cmaplev}
For any \w{\bZ\in\Sa} and \w{R=\Fp} or $\QQ$, the mapping space functor \w{\mapa(\bZ,-)}
is a level homotopy functor \w[.]{\SR\to\Sa}
\end{cor}

\end{document}